\documentclass[reqno,12pt,a4paper]{amsart}

\usepackage{tikz}

\usepackage[mathscr]{eucal}
\usepackage{graphics,epic}
\usepackage{amsfonts}
\usepackage{amscd}
\usepackage{latexsym}
\usepackage{amsmath,amssymb, amsthm, stmaryrd, bm, bbm}
\usepackage[all,2cell]{xy}
\usepackage{rotating}

\newtheorem{theorem}{Theorem}[section]
\newtheorem*{theorem*}{Theorem}

\newtheorem{lemma}[theorem]{Lemma}
\newtheorem{proposition}[theorem]{Proposition}
\newtheorem{corollary}[theorem]{Corollary}

\newtheorem*{conjecture*}{Conjecture}

\newtheorem*{question*}{Question}
\theoremstyle{remark}
\newtheorem{remark}[theorem]{Remark}
\newtheorem{example}[theorem]{Example}
\theoremstyle{definition}

\newcommand{\ie}{{\em i.e.~}\ }

\newcommand{\opname}[1]{\operatorname{\mathsf{#1}}}

\renewcommand{\mod}{\opname{mod}\nolimits}

\newcommand{\grmod}{\opname{grmod}\nolimits}
\newcommand{\proj}{\opname{proj}\nolimits}

\newcommand{\Grmod}{\opname{Grmod}\nolimits}

\newcommand{\grproj}{\opname{grproj}}

\newcommand{\rep}{\opname{rep}\nolimits}
\newcommand{\Rep}{\opname{Rep}\nolimits}

\newcommand{\pretr}{\opname{pretr}\nolimits}

\newcommand{\thick}{\opname{thick}\nolimits}

\newcommand{\per}{\opname{per}\nolimits}

\newcommand{\id}{\mathrm{id}}
\newcommand{\Id}{\mathrm{Id}}

\newcommand{\Hom}{\opname{Hom}}

\newcommand{\cHom}{{\mathcal{H}\it{om}}}

\newcommand{\ten}{\otimes}

\newcommand{\Cone}{\opname{Cone}}

\newcommand{\Tot}{\opname{Tot}}

\newcommand{\ca}{{\mathcal A}}
\newcommand{\cb}{{\mathcal B}}
\newcommand{\cc}{{\mathcal C}}
\newcommand{\cd}{{\mathcal D}}

\newcommand{\ch}{{\mathcal H}}
\newcommand{\ci}{{\mathcal I}}

\newcommand{\cp}{{\mathcal P}}

\newcommand{\ct}{{\mathcal T}}

\newcommand{\cx}{{\mathcal X}}
\newcommand{\cy}{{\mathcal Y}}

\renewcommand{\hat}[1]{\widehat{#1}}

\numberwithin{equation}{section}

\begin{document}

\title[The derived Hall algebra of a graded gentle one-cycle algebra I]{On the derived Hall algebra of a graded gentle one-cycle algebra I: the triangle structure}

\author{Hui Chen, and Dong Yang}

\address{Hui Chen, School of Biomedical Engineering and Informatics, Nanjing Medical University, Nanjing 211166, P. R. China.}
\email{huichen@njmu.edu.cn}

\address{Dong Yang, School of Mathematics, Nanjing University, Nanjing 210093, P. R. China}
\email{yangdong@nju.edu.cn}

\begin{abstract} Under a mild condition, the perfect derived category and the finite-dimensional derived category of a graded gentle one-cycle algebra are described as twisted root categories of certain infinite quivers of type $\mathbb{A}_\infty^\infty$. As a consequence, it is shown that if $\ct$ is the perfect (respectively, finite-dimensional) derived category of such a graded gentle one-cycle algebra, then its triangle structure is up to triangle equivalence determined by the underlying additive category.
\end{abstract}
\maketitle

\section{Introduction}
\label{s:introduction}

This is the first paper of a series of papers which studies the derived Hall algebra of a graded gentle one-cycle algebra. For graded algebras of dual numbers and (ungraded) gentle one-cycle algebras of infinite global dimension a description of the derived Hall algebra can be found in \cite{KellerYangZhou09} and in \cite{BobinskiSchmude20}, respectively.

Gentle algebras originated in the study of derived categories of quivers of type $\mathbb{A}$ and of type $\tilde{\mathbb{A}}$~\cite{AssemHappel81,AssemSkowronski87} and appeared also in cluster theory \cite{AssemBruestleCharbonneauPlamondon10,Labardini09}. They are finite-dimensional algebras with nice homological properties, for example, they are Gorenstein \cite{GeissReiten05} and they are Koszul since they are quadratic monomial. A remarkable result of Vossieck \cite{Vossieck01} states that a finite-dimensional algebra over an algebraically closed field is derived-discrete if and only if it is either derived equivalent to the path algebra of a Dynkin quiver or a gentle one-cycle algebra which does not satisfy the clock condition. Recently, graded gentle algebras appeared in the study of topological Fukaya categories of graded marked surfaces~\cite{HaidenKatzarkovKontsevich17,LekiliPolishchuk20}. From the results in \cite{KalckYang18} graded gentle one-cycle algebras look `generically' `derived-discrete'. For the non-positive ones this is made rigorous in \cite{FushimiXing26} and a certain converse is also given there. 
This paper continues the study of derived categories of graded gentle one-cycle algebras in \cite{KalckYang18}. The main result is the following characterisation.

\begin{theorem}[{Theorem~\ref{thm:perfect-derived-category-as-twisted-root-category}}]
\label{thm:main}
The following conditions are equivalent for a triangulated category $\ct$:
\begin{itemize}
\item[(i)] $\ct$ is the perfect derived category of a graded gentle one-cycle algebra,
\item[(ii)] $\ct$ is triangle equivalent to one triangulated category in the following four families:
\begin{itemize}
\item[(1)] $\cd^b(\rep^b(Q^l))/\Sigma^d\sigma^{-r}$, where $d\geq 1$ and $r\in\mathbb{Z}\backslash\{0\}$, 
\item[(2)] $\cd^b(\rep^b(Q_{p,q}))/\Sigma^d\sigma_{p,q}^{-r}$, where $d,p,q\geq 1$ and $r\in\mathbb{Z}\backslash\{0\}$, 
\item[(3)] $\cd^b(\rep_{\rm nil}^b(Q))$, where $Q$ is a cyclic quiver,
\item[(4)] $\cd^b(\rep^b(Q))$, where $Q$ is a finite quiver of type $\tilde{\mathbb{A}}_{p,q}$, $p,q\geq 1$.
\end{itemize}
\end{itemize}
\end{theorem}
Here in the theorem $Q^l$ and $Q_{p,q}$ are two infinite quivers of type $\mathbb{A}_\infty^\infty$, $\sigma$ and $\sigma_{p,q}$ are quiver-automorphisms (see Sections~\ref{ss:twisted-root-category-linear-orientation}~and~\ref{ss:twisted-root-category-zigzag-orientation} for the precise definitions), $\cd^b$ is the bounded derived category, $\Sigma$ is the shift functor, $\rep^b$ is the category of finite-dimensional representations and $\rep_{\rm nil}^b(Q)$ is the category of finite-dimensional nilpotent representations. In the families (1) and (2), the triangle structures on the orbit categories are given by \cite[Theorem 9.9]{Keller05}. They are known as $d$-root categories when $r=0$, therefore we call them twisted $d$-root categories. For a triangulated category $\ct$ in (ii), the family number with the integers $p,q,d,r$ determines the AG-invariant of the corresponding graded gentle one-cycle algebra introduced in \cite{Avella-AlaminosGeiss08,LekiliPolishchuk20}, and vice versa. There is also a version of Theorem~\ref{thm:main} with the perfect derived category replaced by the finite-dimensional derived category (Theorem~\ref{thm:finite-dimensional-derived-category-as-twisted-root-category}). This characterisation provides an interesting connection between different graded gentle one-cycle algebras (Corollary~\ref{cor:comparing-graded-gentle-one-cycle-algebras}), and also shows that  if $\ct$ is the perfect derived category (respectively, the finite-dimensional derived category) of a graded gentle one-cycle algebra, then its triangle structure is up to triangle equivalence determined by the underlying additive category.

\begin{corollary}[{Corollary~\ref{cor:additive=>triangle}}]
\label{cor:main}
Let $A$ and $A'$ be two graded gentle one-cycle algebras. Then the following conditions are equivalent:
\begin{itemize}
\item[(i)] there is a triangle equivalence $\per(A)\to\per(A')$,
\item[(ii)] there is a $k$-linear equivalence $\per(A)\to\per(A')$,
\item[(iii)] there is a triangle equivalence $\cd_{fd}(A)\to\cd_{fd}(A')$,
\item[(iv)] there is a $k$-linear equivalence $\cd_{fd}(A)\to\cd_{fd}(A')$.
\end{itemize}
\end{corollary}

This paper is structured as follows. In Section~\ref{s:quivers}, we recall some results on the representation theory on infinite quivers and finite graded quivers, and discuss the Hom-finiteness of the perfect derived category of a graded path algebra. In Section~\ref{s:triangulated-orbit-categories}, we recall Keller's results on triangulated orbit categories. In Section~\ref{s:twisted-root-categories}, we introduce twisted root categories and study in detail twisted root categories of $Q^l$ and $Q_{p,q}$. In Section~\ref{s:derived-category-of-graded-path-algebra}, we describe the perfect derived category of the graded path algebra of a graded quiver as a certain twisted root category. In Section~\ref{s:graded-gentle-one-cycle-algebras}, we prove Theorem~\ref{thm:main} and Corollary~\ref{cor:main}.

Throughout, let $k$ be a field and $\mathbb{N}$ denote the set of positive integers.

\medskip
\noindent{\it Acknowledgements.}
The authors are grateful to Bangming Deng and Haicheng Zhang for answering their questions, and to Bernhard Keller for suggesting the terminology ``strictly pretriangulated category" and ``strict dg enhancement". They are indebted to Zhengfang Wang for sharing his knowledge on graded gentle algebras. They thank the referee for carefully reading the manuscript and for his helpful remarks, especially for pointing out an error. The second author acknowledges support by National Key R\&D Program of China 2024YFA1013801 and the National Natural Science Foundation of China No. 12031007. 

\section{Quivers}
\label{s:quivers}
In this section we recall some results on the representation theory of quivers and graded quivers.

\subsection{Quivers}\label{ss:quiver}

We follow \cite{BautistaLiuPaquette13}. Let $Q$ be a (possibly infinite) quiver with vertex set $Q_0$ and arrow set $Q_1$. For an arrow $\alpha$, we denote by $s(\alpha)$ and $t(\alpha)$ the source and the target of $\alpha$, respectively. A \emph{non-trivial path} in $Q$ is a sequence of arrows $\alpha_l\cdots\alpha_1$ with $s(\alpha_{m+1})=t(\alpha_m)$ for all $1\leq m\leq l-1$. An \emph{oriented cycle} is a non-trivial path $\alpha_l\cdots\alpha_1$ with $t(\alpha_l)=s(\alpha_1)$. For each $i\in Q_0$, we have a \emph{trivial path} $e_i$ with $s(e_i)=t(e_i)=i$.

A representation $V$ of $Q$ consists of a collection $(V_i\mid i\in Q_0)$ of vector spaces over $k$ and a collection $(V_\alpha\colon V_{s(\alpha)}\to V_{t(\alpha)}\mid \alpha\in Q_1)$ of $k$-linear maps. It is said to be \emph{locally finite-dimensional} if $V_i$ is finite-dimensional for all $i\in Q_0$ and \emph{finite-dimensional} if $\bigoplus_{i\in Q_0}V_i$ is finite-dimensional.  For example, the simple representation $S_i$ supported at $i\in Q_0$, defined by $(S_i)_i=k$ and $(S_i)_j=0$ for $j\neq i$ and $(S_i)_\alpha=0$ for $\alpha\in Q_1$, is a finite-dimensional representation. 
Let $\Rep(Q)$ denote the category of representations of $Q$, $\rep(Q)$ the subcategory of $\Rep(Q)$ consisting of locally finite-dimensional representations of $Q$ and 
$\rep^b(Q)$ the subcategory of finite-dimensional representations.

\begin{example}
\label{ex:representations-of-A-double-infinity}
Let $Q$ be a quiver of type $\mathbb{A}_\infty^\infty$:
\[\xymatrix{\ldots\ar@{-}[r]&i-1\ar@{-}[r]&i\ar@{-}[r]&i+1\ar@{-}[r]&\ldots}.\]
For $a,b\in\{-\infty\}\cup\mathbb{Z}\cup\{+\infty\}$ with $a\leq b$, we construct a representation $V_{a,b}$ of $Q$ by setting
the $k$-vector space $V_{a,b}(i)$ associated to a vertex $i$ as
\[V_{a,b}(i)=\begin{cases} k & \text{if } a\leq i\leq b\\ 0 & \text{otherwise}\end{cases}\]
and setting the $k$-linear map $V_{a,b}(\alpha)$ associated to an arrow $\alpha\colon i\rightarrow j$ as
\[V_{a,b}(\alpha)=\begin{cases} \mathrm{id}_k &\text{if } a\leq i,j\leq b\\ 0 & \text{otherwise}.\end{cases}\]
These representations are indecomposable with endomorphism algebra isomorphic to $k$ and every indecomposable representation of $Q$ is isomorphic to one such. The simple representations are precisely the $S_i=V_{i,i}$'s ($i\in\mathbb{Z}$) and the finite-dimensional indecomposable representations are precisely the $V_{a,b}$'s with $a,b\in\mathbb{Z}$ and $a\leq b$. See for example \cite[Section 5]{BautistaLiuPaquette13}.
\end{example}

For each $i\in Q_0$, we define a representation $P_i$ of $Q$ as follows. For $j\in Q_0$, the vector space $(P_i)_j$ is the vector space spanned by paths from $i$ to $j$; and for $\alpha\in Q_1$, the $k$-linear map $(P_i)_{\alpha}\colon (P_i)_{s(\alpha)}\to (P_i)_{t(\alpha)}$ takes a path $p$ to $\alpha p$. Let $\proj (Q)$ denote the full subcategory of $\Rep(Q)$ consisting of finite direct sums of $P_i$'s, $i\in Q_0$. A representation $V$ of $Q$ is said to be \emph{finitely presented} if there is a short exact sequence $0\to P^{-1}\to P^0\to V\to 0$, where $P^{-1},~P^0\in \proj(Q)$. Let $\rep^+(Q)$ be the full subcategory of $\Rep(Q)$ consisting of finitely presented representations. When $Q$ is a finite quiver (\ie both $Q_0$ and $Q_1$ are finite sets) , the categories $\rep(Q)$ and $\rep^b(Q)$ coincide. If further $Q$ has no oriented cycles, they also coincide with $\rep^+(Q)$.

$Q$ is said to be \emph{strongly locally finite} if for any vertex $i$, the number of vertices $j$ such that there are arrows between $i$ and $j$ are finite, and for any two vertices $i$ and $j$, the number of paths from $i$ to $j$ is finite. It is clear that $\rep^b(Q)$ is a Hom-finite and Ext-finite hereditary abelian $k$-category.

\begin{proposition}[{\cite[Proposition 1.15]{BautistaLiuPaquette13}}]
\label{prop:Hom-finiteness-of-finitely-presented-representations}
Assume that $Q$ is strongly locally finite. Then $\rep^+(Q)$ is a Hom-finite and Ext-finite hereditary abelian $k$-category.
\end{proposition}

If $\sigma$ is a quiver-automorphism of $Q$, then the corresponding push-out functor $\sigma_*$ is a $k$-linear automorphism of $\Rep(Q)$, which restricts to $k$-linear automorphisms of $\rep(Q)$, $\rep^+(Q)$ and $\rep^b(Q)$. Precisely, for $V\in\Rep(Q)$,  we have $(\sigma_*(V))_i=V_{\sigma^{-1}(i)}$ for $i\in Q_0$ and $(\sigma_*(V))_\alpha=V_{\sigma^{-1}(\alpha)}$. For example $\sigma_*(S_i)=S_{\sigma(i)}$ and $\sigma_*(P_i)=P_{\sigma(i)}$. By abuse of notation, we will often denote $\sigma_*$ and its derived functor also by $\sigma$.

\subsection{Graded quivers}
\label{ss:graded-quiver}
Most of this subsection is taken from \cite[Section 4]{Yang18}. We will state the results and use them for the quiver which is `Koszul dual' to the quiver in \cite[Section 4]{Yang18}. Thus in order to avoid confusion we include them in detail.

Let $Q$ be a finite graded quiver, \ie a finite quiver together with a map which assigns to each arrow $\alpha$ an integer $|\alpha|$ called its degree. The degree of a trivial path is $0$ and the degree of a non-trivial path $\alpha_l\cdots\alpha_1$ is $|\alpha_1|+\ldots+|\alpha_l|$. We consider an (ungraded) quiver as a graded quiver whose arrows are in degree $0$.  The \emph{graded path algebra} $kQ$ of $Q$ is the graded $k$-algebra which has all paths in $Q$ (including the trivial paths) as basis and the multiplication is given by concatenation of paths. So the degree $m$ component of $kQ$ has all paths of degree $m$ as basis. The \emph{complete path algebra} $\widehat{kQ}$ of $Q$ is the completion of $kQ$ in the category of graded $k$-algebras with respect to the $J$-adic topology, where $J$ is the ideal of $kQ$ generated by all arrows. So the degree $m$ component of $\widehat{kQ}$ consists of (possibly infinite) linear combinations of all paths of degree $m$. For example, if $Q$ is the Jordan quiver with the unique arrow in degree $-1$, then $\hat{kQ}=kQ$ is isomorphic to the graded polynomial algebra $k[x]$ with $x$ in degree $-1$. We refer to \cite[Section II.3]{CaenepeelVanOystaeyen88} for the theory on completions of graded rings.

For $\alpha\in Q_1$, we introduce its formal inverse $\alpha^{-1}$, namely, $s(\alpha^{-1})=t(\alpha)$,  $t(\alpha^{-1})=s(\alpha)$, $|\alpha^{-1}|=-|\alpha|$ and $\alpha\alpha^{-1}=e_{t(\alpha)}$, $\alpha^{-1}\alpha=e_{s(\alpha)}$.  
A \emph{walk} is a formal product $w=c_l\cdots c_1$, where each $c_m$ is either a trivial path, or an arrow, or the inverse of an arrow such that $s(c_{m+1})=t(c_m)$ for $1\leq m\leq l-1$. Define $s(w)=s(c_1)$ as the source of $w$ and $t(w)=t(c_l)$ as the target of $w$. We say that $w$ is a \emph{closed walk} if $s(w)=t(w)$. The \emph{degree} of $w$ is defined as $|w|=|c_l|+\ldots+|c_1|$. The \emph{cycle degree} of $Q$ is defined as $d=0$ if all closed walks are of degree $0$, and otherwise as the minimal positive integer $d$ such that there are closed walks of degree $d$.   Note that the degree of any closed walk is a multiple of $d$, if $Q$ is connected. For example, the cycle degree of the Jordan quiver with the unique arrow in degree $-1$ is  $1$.

\medskip
Let $\Grmod(kQ)$ be the category of graded $kQ$-modules, $\grmod(kQ)$ the full subcategory of finitely presented graded modules and $\grmod_0(kQ)$ the full subcategory of finite-dimensional graded modules. $\Grmod(kQ)$ is a hereditary abelian $k$-category, and  by \cite[Section 1.4]{Smith12}, $\grmod(kQ)$ is a hereditary abelian $k$-category.  It follows that $\grmod_0(kQ)$ is also a hereditary abelian $k$-category; moreover, it is Hom-finite and Ext-finite. We want to study these categories in terms of quiver representations. For this purpose an infinite (ungraded) quiver $R$ together with a quiver-automorphism $\sigma$ is defined. The vertices of $R$ are $(i,j)$, where $i\in Q_0$ and $j\in\mathbb{Z}$. The arrows of $R$ are of the form $(\alpha,j)\colon (t(\alpha),j)\to (s(\alpha),j+|\alpha|)$, where $\alpha\in Q_1$ and $j\in\mathbb{Z}$. Then the assignments $(i,j)\mapsto (i,j-1)$ and $(\alpha,j)\mapsto(\alpha,j-1)$ defines a quiver-automorphism $\sigma$ of $R$.  The $\sigma$ here is the inverse of the $\sigma$ in \cite[Section 4]{Yang18} (we make this change in order to use the push-out functor $\sigma_*$ instead of the pull-back functor $\sigma^*$).

\begin{example}\label{ex:kronecker-1}
Let $Q$ be the graded Kronecker quiver
\[
\xymatrix@C=3pc{
1& 2\ar@<.7ex>[l]^\alpha \ar@<-.7ex>[l]_\beta 
}
\]
with $|\alpha|=0$ and $|\beta|=2$. Then the quiver $R$ is
\[{\scriptsize
\xymatrix@R=1pc@C=3pc{
\ar[ddrr] & (1,-2)\ar[dd]|(0.65){(\alpha,-2)} \ar[ddrr]|(0.3){(\beta,-2)} & (1,-1)\ar[dd]|(0.65){(\alpha,-1)} \ar[ddrr]|(0.3){(\beta,-1)} & (1,0)\ar[dd]|(0.65){(\alpha,0)} \ar[ddrr]|(0.3){(\beta,0)} & (1,1)\ar[dd]|(0.65){(\alpha,1)} \ar[ddrr]|(0.3){(\beta,1)}  & (1,2)\ar[dd]|(0.65){(\alpha,2)}  & \\
\ldots & & & & & & \ldots\\
& (2,-2) & (2,-1) & (2,0) & (2,1) & (2,2) &  &&
}
}
\]
and $\sigma$ is the shift to the left by one step.
\end{example}

\begin{lemma}
\label{lem:from-graded-modules-to-representations}
There is a $k$-linear equivalence $F\colon\Grmod(kQ)\to\Rep(R)$ such that $F\circ\langle 1\rangle=\sigma_*\circ F$. This equivalence restricts to $k$-linear equivalences $\grmod(kQ)\to \rep^+(R)$ and $\grmod_0(kQ)\to \rep^b(R)$.
\end{lemma}
\begin{proof}
For $M\in \Grmod(kQ)$, define $F(M)$ as the representation $V$ of $R$ with $V_{(i,j)}=M^j e_i$ and $V_{(\alpha,j)}\colon V_{(t(\alpha),j)}\to V_{(s(\alpha),j+|\alpha|)}$ being given by right action by $\alpha$. It is easy to check that $F$ is a $k$-linear equivalence and $F(M\langle 1\rangle)=\sigma_*(F(M))$. 
That $F$ restricts to a $k$-linear equivalence $\grmod_0(kQ)\to \rep^b(R)$ is clear. Finally, observe that $F(e_i kQ\langle j\rangle)= P_{(i,-j)}$ and hence $F$ restricts to a $k$-linear equivalence $\grproj(kQ)\to\proj(R)$, and hence to a $k$-linear equivalence $\grmod(kQ)\to\rep^+(R)$.
\end{proof}

Assume that $Q$ is connected. We define a (possibly infinite) (ungraded) quiver $\tilde{Q}$. Fix a vertex $i$ of $Q$. The vertices of $\tilde{Q}$ are equivalence classes of walks of $Q$ whose target is $i$, where two walks $u$ and $v$ are equivalent if they have the same source and are of the same degree. For a walk $u$ with target $i$ we denote by $[u]$ the equivalence class of $u$. The arrows of $\tilde{Q}$ are of the form $\alpha^{[u]}\colon [u]\to [u\alpha]$, where $[u]$ runs over all vertices of $\tilde{Q}$ and $\alpha$ runs over all arrows of $Q$ such that $t(\alpha)=s(u)$. Let $w$ be a closed walk with source $i$ and of degree $-d$. Then the assignments $[u]\mapsto [wu]$ and $\alpha^{[u]}\mapsto \alpha^{[wu]}$ defines a quiver-automorphism $\mathrm{s}$ of $\tilde{Q}$.
This $\mathrm{s}$ is the inverse of the $s$ in \cite[Section 4]{Yang18}.

\begin{example}\label{ex:kronecker-2}
Let $Q$ be the graded Kronecker quiver
\[
\xymatrix@C=3pc{
1& 2\ar@<.7ex>[l]^\alpha \ar@<-.7ex>[l]_\beta 
}
\]
with $|\alpha|=0$ and $|\beta|=2$. Take $i=1$. Then $d=2$ and $\tilde{Q}$ is the quiver
\[{\scriptsize
\xymatrix@R=1pc@C=3pc{
\ar[ddr] & [\alpha\beta^{-1}]\ar[dd]|{\alpha^{[\alpha\beta^{-1}]}}\ar[ddr]|{\beta^{[\alpha\beta^{-1}]}} & [e_1]\ar[dd]|{\alpha^{[e_1]}}\ar[ddr]|{\beta^{[e_1]}} & [\beta\alpha^{-1}] \ar[dd]|{\alpha^{[\beta\alpha^{-1}]}}\ar[ddr]|{\beta^{[\beta\alpha^{-1}]}}& [\beta\alpha^{-1}\beta\alpha^{-1}]\ar[dd]|{\alpha^{[\beta\alpha^{-1}\beta\alpha^{-1}]}} \ar[ddr]& & \\
 \ldots&&&&&\ldots\\
 & [\alpha\beta^{-1}\alpha] & [\alpha] & [\beta] & [\beta\alpha^{-1}\beta] & 
}
}
\]
and $\mathrm{s}$ is the shift to the left by one step.
%
\end{example}

The following lemma is \cite[Lemma 4.10]{Yang18}. 

\begin{lemma}
\label{lem:covering-quiver}
If $d=0$, then $\tilde{Q}$ is isomorphic to the underlying quiver of $Q^{op}$, $\mathrm{s}$ is the identity, and $R$ is isomorphic to the disjoint union of $\mathbb{Z}$ copies of $\tilde{Q}$. If $d\geq 1$, then $R$ is isomorphic to the disjoint union of $d$ copies of $\tilde{Q}$.
\end{lemma}

Let us describe the isomorphism when $d\geq 1$. The disjoint union $\coprod_d \tilde{Q}$ of $d$ copies of $\tilde{Q}$ has vertices $([u],j)$, where $[u]$ runs over all equivalence classes of walks of $Q$ with target $i$ and $0\leq j\leq d-1$, and arrows $(\alpha^{[u]},j)$, where $[u]$ runs over all equivalence classes of walks of $Q$, $\alpha$ runs over all arrows of $Q$ with $t(\alpha)=s(u)$, and $0\leq j\leq d-1$. Then the assignments $([u],j)\mapsto (s(u),j+|u|)$ and $(\alpha^{[u]},j)\mapsto (\alpha,j+|u|)$ defines a quiver-isomorphism $\varphi\colon\coprod_d \tilde{Q}\to R$. 
The quiver-automorphism $\varphi^{-1}\sigma\varphi$ of $\coprod_d \tilde{Q}$, still denoted by $\sigma$ by abuse of notation, takes $([u],j)$ to $([u],j-1)$ and $(\alpha^{[u]},j)$ to $(\alpha^{[u]},j-1)$ for $1\leq j\leq d-1$, and takes $([u],0)$ to $(\mathrm{s}([u]),d-1)$ and $(\alpha^{[u]},0)$ to $(\mathrm{s}(\alpha^{[u]}),d-1)$.

\subsection{Hom-finiteness for the perfect derived category of a graded path algebra}

First we quickly introduce dg algebras and their derived categories. For details we refer to \cite{Keller94,Keller06d,KalckYang18,KalckYang18a}. 

A dg $k$-algebra $A$  is a graded $k$-algebra $A=\bigoplus_{i\in\mathbb{Z}}A^i$ endowed with a differential $d_A$ of degree $1$ such that the graded Leibniz rule holds
\[d_A(aa')=d_A(a)a'+(-1)^{|a|}ad_A(a'),\]
where $a$ is homogeneous of degree $|a|$. For example, a graded $k$-algebra can be viewed as a dg $k$-algebra with trivial differential. 
A (right) \emph{dg
$A$-module} is a (cochain) complex $M$ (of $k$-vector spaces) endowed with an
$A$-action from the right
\[M^i\ten A^j\rightarrow M^{i+j}, m\ten a\mapsto ma\]
such that the graded Leibniz rule holds
\[d_M(ma)=d_M(m)a+(-1)^{|m|}md_A(a),\]
where $m\in M$ is homogeneous of degree $|m|$.

For a dg $k$-algebra $A$, we denote by $\cd(A)$ the derived category of dg $A$-modules. It is a triangulated $k$-category whose objects are dg $A$-modules and whose shift functor is the shift of complexes. Morphism spaces in $\cd(A)$ can be computed by
\[
\Hom_{\cd(A)}(M,N)\cong H^0\cHom_A(\mathbf{p}M,N),
\]
where $\mathbf{p}M$ is a $K$-projective dg $A$-module which is quasi-isomorphic to $M$ ($A$ itself is a typical $K$-projective dg $A$-module), and for two dg $A$-modules $M$ and $N$, $\cHom_A(M,N)$ is the complex whose degree $n$ component $\cHom_A^n(M,N)$ consists of homogeneous homomorphisms $f\colon M\to N$ of graded $A$-modules of degree $n$ and the differential takes such an $f$ to $d_N\circ f-(-1)^n f\circ d_M$.

Let $\per(A)$ denote the perfect derived category, \ie the smallest triangulated subcategory of $\cd(A)$ which contains $A$ and which is closed under taking direct summands, and $\cd_{fd}(A)$ denote the finite-dimensional derived category, \ie $\cd_{fd}(A)$ is the full subcategory of $\cd(A)$ which consists of dg modules $M$ with finite-dimensional total cohomology.
A triangle equivalence $\cd(B)\rightarrow\cd(A)$ between derived categories of dg algebras 
restricts to triangle equivalences $\per(B)\rightarrow\per(A)$ and $\cd_{fd}(B)\rightarrow\cd_{fd}(A)$.

Next let $Q$ be a finite graded quiver and view its graded path algebra $kQ$ as a dg algebra with trivial differential. Since $\per(kQ)$ is by definition idempotent-complete, the condition (iv) below is equivalent to: $\per(kQ)$ is Hom-finite and Krull--Schmidt.

\begin{lemma}\label{lem:finiteness-of-the-perfect-derived-category}
The following conditions are equivalent:
\begin{itemize}
\item[(i)] $\widehat{kQ}=kQ$,
\item[(ii)] for any $p\in\mathbb{Z}$, $Q$ has only finitely many paths of degree $p$,
\item[(iii)] for any $p\in\mathbb{Z}$, the degree $p$ component of $kQ$ is finite-dimensional over $k$,
\item[(iv)] $\per(kQ)$ is Hom-finite,
\item[(v)] $\widetilde{Q}$ is strongly locally finite,
\item[(vi)] $R$ is strongly locally finite.
\end{itemize}
\end{lemma}
\begin{proof}
The equivalences (i)$\Leftrightarrow$(ii)$\Leftrightarrow$(iii) are clear. The equivalence (iii)$\Leftrightarrow$(iv) follows by d\'evissage from the equalities
\[
\Hom_{\per(kQ)}(kQ,\Sigma^ p kQ)=H^p\cHom_{kQ}(kQ,kQ)=H^p(kQ)=(kQ)^p,
\]
since $kQ$ is a classical generator of $\per(kQ)$. For the equivalence (ii)$\Leftrightarrow$(vi), it is enough observe that the paths (respectively, arrows) in $R$ from a vertex $(i,j)$ to another vertex $(i',j')$ are in bijection with the paths (respectively, arrows) in $Q$ from $i$ to $i'$ of degree $j'-j$. Finally, the equivalence (v)$\Leftrightarrow$(vi) follows from Lemma~\ref{lem:covering-quiver}.
\end{proof}

It is clear that the condition (ii) in Lemma~\ref{lem:finiteness-of-the-perfect-derived-category} implies that $Q$ has no oriented cycles of degree $0$. The converse is not true. For example, let $Q$ be the graded quiver
\[\xymatrix{1\ar@(ul,dl)_{t_1}\ar[r]^{\alpha} & 2\ar@(ur,dr)^{t_2}}\]
with $|t_1|=-3$, $|\alpha|=1$ and $|t_2|=3$. Then $Q$ has no oriented cycles of degree $0$, in fact, there are only $2$ paths of degree $0$ (the two trivial paths), but there are infinitely many paths of degree $1$: $t_2^n \alpha t_1^n$, $n\in\mathbb{N}$.

\section{Triangulated orbit categories}
\label{s:triangulated-orbit-categories}

In this section, we recall a special case of Keller's result on triangulated orbit categories in \cite[Section 9]{Keller05}, and make it explicit that a triangle functor between triangulated categories with autoequivalences which admits a `dg lift' induces a triangle functor between the triangulated orbit categories.

\subsection{Pretriangulated dg categories}
\label{ss:dg-enhancements}

Let $\ca$ be a \emph{dg $k$-category}, \ie a $k$-category in which all
the $k$-vector spaces $\Hom_{\ca}(X,Y)$ are endowed with the
structure of a complex of $k$-vector spaces such that the compositions
\[\Hom_{\ca}(Y,Z)\ten_k \Hom_{\ca}(X,Y)\rightarrow \Hom_{\ca}(X,Z)\]
are chain maps of complexes. We define $H^0\ca$ to be the
$k$-category whose objects are the same as $\ca$ and whose morphism
space $\Hom_{H^0\ca}(X,Y)$ is the 0-th cohomology of the complex
$\Hom_{\ca}(X,Y)$. A \emph{dg functor} $F\colon\ca\rightarrow\cb$ between
dg $k$-categories is a $k$-linear functor such that each map
$F(X,Y)\colon\Hom_{\ca}(X,Y)\rightarrow\Hom_{\cb}(FX,FY)$ is a chain map
of complexes. We define $H^0F$ to be the $k$-linear functor
$H^0F\colon H^0\ca\rightarrow H^0\cb$ such that $H^0FX=FX$ and
$(H^0F)(X,Y)=H^0(F(X,Y))$. 
A dg functor $F\colon \ca\to\cb$ is called an \emph{equivalence} if it is fully faithful and dense, and a \emph{quasi-equivalence} if $H^0F\colon H^0\ca\to H^0\cb$ is a $k$-linear equivalence. Note that this definition of quasi-equivalence is different from the one in \cite[Section 7]{Keller94}.

Below we list a few facts on pretriangulated dg categories, for details see \cite{BondalKapranov90,Drinfeld04,Keller05,Keller06d,KalckYang18}. Here what we call `strictly pretriangulated dg category' is called `+-pretriangulated dg category' in \cite{BondalKapranov90}, `exact dg category' in \cite{Keller99}, `strongly pretriangulated dg category' in \cite{Drinfeld04}, and `pretriangulated dg category' in \cite{Keller05,Keller06d,KalckYang18}.

\begin{itemize}
\item
A \emph{strictly pretriangulated dg $k$-category} (respectively, \emph{pretriangulated dg $k$-category}) is roughly speaking a dg $k$-category which has shifts and cones (respectively, up to homotopy). $H^0\ca$ is naturally a triangulated category.
\item
For a strictly pretriangulated dg $k$-category $\ca$, taking an object to its shift extends to a dg automorphism of $\ca$, which we denote by $[1]$. And $H^0\ca$ is a triangulated $k$-category with shift functor 
$\Sigma=H^0[1]$ and with triangles the $\Sigma$-sequences isomorphic to
$X\stackrel{f}{\rightarrow} Y\rightarrow \Cone(f)\rightarrow X[1]$, 
where $f$ runs over all closed morphisms of degree $0$ in $\ca$.
\item 
Let $\ca$ and $\cb$ be two strictly pretriangulated dg $k$-categories and $F\colon\ca\to\cb$ be a dg functor. Then there is a canonical  isomorphism $\varphi_F\colon F[1]\to [1]F$ of dg functors and $F$ also commutes with cones. As a consequence, $(H^0F,H^0\varphi_F)\colon H^0(\ca)\to H^0(\cb)$ is a triangle functor. We will usually omit $H^0\varphi_F$.
\item 
Let $\ca$ be a strictly pretriangulated dg $k$-category and $\ca'$ be a full dg subcategory. The \emph{pretriangulated hull} $\pretr(\ca')$ of $\ca'$ is the smallest pretriangulated dg subcategory of $\ca$ which contains $\ca'$. The embedding $\ca'\to\pretr(\ca')$ is universal among all dg functors from $\ca'$ to strictly pretriangulated dg $k$-categories, and hence does not depend on $\ca$. If $\ca'$ is pretriangulated, then $\ca'\to\pretr(\ca')$ is a quasi-equivalence.
\item
Let $\ca$ be a pretriangulated dg $k$-category and $\ca'$ be a full dg subcategory. Then the dg quotient $\ca/\ca'$ is again pretriangulated and there is a canonical triangle equivalence $H^0\ca/\thick(H^0\ca')\to H^0(\ca/\ca')$. However, even if $\ca$ is strictly pretriangulated, $\ca/\ca'$ is in general not strictly pretriangulated.
\item 
Let $\ca$ be an additive $k$-category. 
We define
$\cc_{dg}^b(\ca)$ to be the category whose objects are bounded complexes of
objects in $\ca$
and whose morphism complexes are
$\Hom_{\cc_{dg}^b(\ca)}(X,Y)=\cHom_{\ca}(X,Y)$, the complex whose degree $n$ component is
\[\cHom_{\ca}^n(X,Y)=\prod_{m\in\mathbb{Z}}\Hom_{\ca}(X^m,Y^{m+n})\]
and whose differential is given by
$d(f)=d_Y\circ f-(-1)^n f\circ d_X$
for $f\in\cHom_\ca^n(X,Y)$. 
Then the category $\cc_{dg}^b(\ca)$ is a strictly pretriangulated dg $k$-category
and $H^0\cc_{dg}^b(\ca)=\ch^b(\ca)$, the bounded homotopy category of $\ca$.

\end{itemize}

\subsection{Triangulated orbit categories}
\label{ss:triangulated-orbit-categories}
We follow \cite[Section 9]{Keller05} (we remark that \cite[Section 9]{Keller05} is much more general). 
Let $\ca$ be a $k$-category (respectively, a dg $k$-category) and $\Phi$ be a $k$-linear autoequivalence (respectively, a dg autoequivalence) of $\ca$. 
First assume that  there exists an isomorphism $\epsilon\colon \Phi^d\to \Id_\ca$ for some $d\in\mathbb{N}$, which is compatible with $\Phi$, that is, $\epsilon\Phi=\Phi\epsilon$. We introduce the \emph{$d$-orbit category} $\ca/_{\hspace{-3pt}d}\Phi$ as the $k$-category (respectively, dg $k$-category) which has the same objects as $\ca$ and for two objects $X$ and $Y$ the morphism space (respectively, complex) is
\[
\Hom_{\ca/_{\hspace{-2pt}d}\Phi}(X,Y)=\bigoplus_{p=0}^{d-1}\Hom_\ca(X,\Phi^p(Y)).
\]
For two morphisms $f\in\Hom_\ca(X,\Phi^p(Y))$ and $g\in\Hom_\ca(Y,\Phi^q(Z))$, considered as morphisms in $\ca/_{\hspace{-3pt}d}\Phi$, the composition $g\circ f$ is defined as $\Phi^p(g)f$ if $p+q<d$ and as $\Phi^{p+q-d}(\epsilon_Z)\Phi^p(g)f$ if $p+q\geq d$. The canonical projection functor $\pi\colon\ca\to\ca/_{\hspace{-3pt}d}\Phi$, which takes an object $X$ to $X$ and a morphism $f$ to $f$, is a $k$-linear functor (respectively, a dg functor).   

Now in general, the \emph{orbit category} $\ca/\Phi$ by definition has the same objects as $\ca$ and for two objects $X$ and $Y$ the morphism space (respectively, complex) is
\[
\Hom_{\ca/\Phi}(X,Y)=\bigoplus_{p\in\mathbb{Z}}\Hom_\ca(X,\Phi^p(Y)).
\]
For two morphisms $f\in\Hom_\ca(X,\Phi^p(Y))$ and $g\in\Hom_\ca(Y,\Phi^q(Z))$, considered as morphisms in $\ca/\Phi$, the composition $g\circ f$ is defined as $\Phi^p(g)f$. The canonical projection functor $\pi\colon\ca\to\ca/\Phi$, which takes an object $X$ to $X$ and a morphism $f$ to $f$, is a $k$-linear functor (respectively, a dg functor).  Here in order to make the composition of morphisms work properly one should fix a quasi-inverse of $\Phi$ and an adjunction between them. We omit the details.

If $\cb$ is another $k$-category (respectively, dg $k$-category), $\Psi$ is a $k$-linear autoequivalence (respectively, dg autoequivalence) of $\cb$, and $F\colon\ca\to\cb$ is a $k$-linear functor (respectively, dg functor) with a fixed isomorphism $\gamma\colon \Psi F\to F\Phi$, then there is an induced $k$-linear functor (respectively, dg functor) $\bar{F}\colon\ca/\Phi\to\cb/\Psi$ such that the following diagram is commutative\footnote{Warning: In the dg setting in the case $(\ca,\Phi)=(\cb,\Psi)$ with $\ca$ being strictly pretriangulated and $F=[1]$, the induced functor $\bar{F}\colon\ca/\Phi\to\ca/\Phi$ may not be the shift functor of $\ca/\Phi$. A misunderstanding of this fact may cause much confusion. In order to obtain the shift functor of $\ca/\Phi$ we have to use the isomorphism $\varphi_\Phi\colon \Phi[1]\to[1]\Phi$.}
\[
\xymatrix{
\ca\ar[r]^F\ar[d]^\pi & \cb\ar[d]^\pi\\
\ca/\Phi \ar[r]^{\bar{F}} & \cb/\Psi.
}
\]
In particular, for $d\in\mathbb{N}$, the $k$-linear autoequivalence (respectively, dg autoequivalence) $\Phi$ of $\ca$, by the trivial commutativity $\Phi\circ\Phi^d=\Phi^d\circ\Phi$, induces a $k$-linear autoequivalence (respectively, a dg autoequivalence) $\bar{\Phi}$ of $\ca/\Phi$, such that there is a natural isomorphism $\epsilon\colon\bar{\Phi}^d\cong \Id_{\ca/\Phi}$ with $\epsilon_X=\id_{\Phi^d(X)}$, which is compatible with $\bar{\Phi}$. It is easy to see that $(\ca/\Phi^d)/_{\hspace{-3pt}d}\bar{\Phi}=\ca/\Phi$ and that there is the following commutative diagram
\begin{align}
\label{cd:finite-orbit-category}
\xymatrix{
\ca\ar[r]^\pi \ar[d]^\pi & \ca/\Phi^d\ar[d]^\pi\\
\ca/\Phi\ar@{=}[r] & (\ca/\Phi^d)/_{\hspace{-3pt}d}\bar{\Phi}.
}
\end{align}

Let $\ct$ be a triangulated $k$-category with a triangle autoequivalence $F$. A \emph{strict dg enhancement}\footnote{As pointed out by Bernhard Keller, a \emph{dg enhancement} of $(\ct,F)$ should be a quadruple $(\ca,\eta,\alpha,\Phi)$, where $\ca$ is a pretriangulated dg $k$-category, $\Phi\colon \ca\to\ca$ is a dg functor, and $(\eta,\alpha)\colon (H^0\ca,H^0\Phi)\to(\ct,F)$ is an equivalence in the 2-category of triangulated categories with autoequivalence.} of $(\ct,F)$ is a triple $(\ca,\eta,\Phi)$, where $\ca$ is a pretriangulated dg $k$-category, $\eta\colon H^0\ca\to\ct$ is a triangle equivalence and $\Phi$ is a dg autoequivalence of $\ca$ such that there is a commutative diagram of triangle functors
\[
\xymatrix{
H^0\ca\ar[r]^{H^0\Phi}\ar[d]^\eta & H^0\ca\ar[d]^\eta\\
\ct\ar[r]^F & \ct.
}
\]
We say that two strict dg enhancements $(\ca,\eta,\Phi)$ and $(\ca',\eta',\Phi')$ of $(\ct,F)$ are \emph{equivalent} if there exists a dg functor $\xi\colon\ca\to\ca'$ such that there are diagrams of dg functors and of triangle functors which are commutative up to isomorphism
\[
\xymatrix{
\ca'\ar[r]^{\Phi'} & \ca'\\
\ca\ar[u]^\xi \ar[r]^{\Phi} & \ca\ar[u]^\xi
}\hspace{25pt}\text{and}\hspace{25pt}
\xymatrix@C=2pc@R=0.7pc{
&H^0\ca'\ar[rr]|(0.4){H^0\Phi'} \ar@{-->}[lddd]|(0.5){\eta'} && H^0\ca'\ar[lddd]|(0.5){\eta'}\\
H^0\ca\ar[rr]|(0.6){H^0\Phi} \ar[dd]|(0.4)\eta \ar[ur]|{H^0\xi} && H^0\ca\ar[dd]|(0.4)\eta \ar[ur]|{H^0\xi}\\
\\
\ct\ar[rr]|(0.6)F && \ct
}
\]

Assume that $(\ct,F)$ admits a strict dg enhancement $(\ca,\eta,\Phi)$. We consider on one hand the dg orbit category $\ca/\Phi$ and on the other hand the orbit category $\ct/F$, which is equivalent to $H^0\ca/H^0\Phi$. It is clear that $H^0\ca/H^0\Phi=H^0(\ca/\Phi)$. Taking $H^0$ of the following two dg functors
\[
\ca\to\ca/\Phi\to\pretr(\ca/\Phi),
\]
we obtain a commutative diagram of $k$-linear functors
\begin{align}
\label{eq:triangulated-hull}
\xymatrix{
H^0\ca\ar[r]^(0.4)\pi \ar[d]_{\eta} & H^0\ca/H^0\Phi\ar[r] \ar[d]_{\eta} & H^0(\pretr(\ca/\Phi))\\
\ct  \ar[r]^(0.4)\pi & \ct/F  \ar[r] & H^0(\pretr(\ca/\Phi)), \ar@{=}[u]
}
\end{align}
the outer square of which is a commutative square of triangle functors. The second functor in the second row is called a \emph{triangulated hull} of $\ct/F$.  When it is a $k$-linear equivalence, it induces a triangle structure on $\ct/F$ such that the projection functor $\pi\colon\ct \to \ct/F$ is a triangle functor. The following result is a special case of \cite[Theorem 9.9]{Keller05}. For an abelian $k$-category $\ch$, let $\cd^b(\ch)_{dg}$ be the standard dg enhancement of $\cd^b(\ch)$, \ie the dg quotient of $\cc^b_{dg}(\ch)$ by its full dg subcategory of acyclic complexes, see \cite[Sections 9.8 and 9.9]{Keller05}. 

\begin{theorem}
\label{thm:triangulated-orbit-category-Keller}
Let $\ch$ be a small Hom-finite and Ext-finite hereditary abelian $k$-category and $F$ be a triangle autoequivalence of $\cd^b(\ch)$. Assume that $F$ satisfies the following conditions:
\begin{itemize}
\item[2)] For each indecomposable $U$ of $\ch$, only finitely many objects $F^i (U)$, $i\in\mathbb{Z}$, lie in $\ch$.
\item[3)] There is an integer $N\geq 0$ such that the $F$-orbit of each indecomposable object of $\cd^b(\ch)$ contains an object $\Sigma^n U$, for some $0\leq n\leq N$ and some indecomposable object $U$ of $\ch$,
\end{itemize}
and that  there exists a dg autoequivalence $\Phi$ of $\cd^b(\ch)_{dg}$ such that $H^0\Phi=F$. 
Then the triangulated hull above is a $k$-linear equivalence. Consequently, the orbit category $\cd^b(\ch)/F$ admits a triangle structure such that the projection functor $\cd^b(\ch)\to\cd^b(\ch)/F$ is a triangle functor. We call this triangle structure the \emph{standard triangle structure} on $\cd^b(\ch)/F$.

\end{theorem}

\subsection{Comparing orbit categories}
\label{ss:comparing-orbit-categories}

Let $\ct$ and $\ct'$ be triangulated $k$-categories with triangle autoequivalences $F$ and $F'$, respectively, and $G\colon \ct\to\ct'$ be a triangle  functor with a fixed isomorphism $GF\cong F'G$. Assume that $(\ct,F)$ and $(\ct',F')$ admit strict dg enhancements $(\ca,\eta,\Phi)$ and $(\ca',\eta',\Phi')$, and that there is a dg functor $\xi\colon\ca\to\ca'$ such that there are diagrams of dg functors  and of triangle functors which are commutative up to isomorphism
\[
\xymatrix{
\ca'\ar[r]^{\Phi'} & \ca'\\
\ca\ar[u]^\xi \ar[r]^{\Phi} & \ca\ar[u]^\xi
}\hspace{25pt}\text{and}\hspace{25pt}
\xymatrix@C=2pc@R=0.7pc{
&H^0\ca'\ar[rr]|(0.4){H^0\Phi'} \ar@{-->}[dd]|(0.4){\eta'} && H^0\ca'\ar[dd]|(0.4){\eta'}\\
H^0\ca\ar[rr]|(0.6){H^0\Phi} \ar[dd]|(0.4)\eta \ar[ur]|{H^0\xi} && H^0\ca\ar[dd]|(0.4)\eta \ar[ur]|{H^0\xi}\\
&\ct'\ar@{-->}[rr]|(0.4){F'} && \ct'\\
\ct\ar[rr]|(0.6)F \ar@{-->}[ur]|{G} && \ct\ar[ur]|{G}
}
\]
The first diagram induces a dg functor $\ca/\Phi\to\ca'/\Phi'$  and further a dg functor $\pretr(\ca/\Phi)\to\pretr(\ca'/\Phi')$. The second diagram further induces a commutative diagram
\[
\xymatrix@C=1.3pc@R=0.9pc{
&H^0\ca'\ar[rr] \ar@{-->}[dd]|(0.4){\eta'} && H^0\ca'/H^0\Phi'\ar@{-->}[dd]|(0.4){\bar{\eta}'}\ar[rr] && H^0(\pretr(\ca'/\Phi'))\ar@{=}[dd]\\
H^0\ca\ar[rr] \ar[dd]|(0.4)\eta \ar[ur]|{H^0\xi} && H^0\ca/H^0\Phi\ar[dd]|(0.4){\bar{\eta}} \ar[ur]|{\overline{H^0\xi}}\ar[rr] && H^0(\pretr(\ca/\Phi))\ar[ur]\ar@{=}[dd]\\
&\ct'\ar@{-->}[rr] && \ct'/F' \ar@{-->}[rr] && H^0(\pretr(\ca'/\Phi'))\\
\ct\ar[rr] \ar@{-->}[ur]|{G} && \ct/F\ar@{-->}[ur]|{\bar{G}} \ar[rr] && H^0(\pretr(\ca/\Phi))\ar[ur]
}
\]
\begin{proposition}
\label{prop:comparing-orbit-categories}
If both triangulated hulls of $\ct/F$ and $\ct'/F'$ above are $k$-linear equivalences, then $\bar{G}\colon\ct/F\to\ct'/F'$ is a triangle functor. Moreover, if $G$ is a triangle equivalence, then $\bar{G}$ is a triangle equivalence. As a consequence, equivalent strict dg enhancements of $(\ct,F)$ induce the same triangle structure on $\ct/F$.
\end{proposition}

\begin{remark}
\label{rem:triangle-structure-on-orbit-category-with-projectives}
Let $\ch$ be a small Hom-finite and Ext-finite hereditary abelian $k$-category with enough projectives and let $\proj\ch$ be the full subcategory of $\ch$ of projective objects. Then $\ca=\cc_{dg}^b(\proj\ch)$ is a strictly pretriangulated dg $k$-category, and taking $H^0$ of the composition $\cc_{dg}^b(\proj\ch)\to\cc_{dg}^b(\ch)\to\cd^b(\ch)_{dg}$ we obtain the canonical functor $K^b(\proj\ch)\to\cd^b(\ch)$, which is a triangle equivalence. Thus by Proposition~\ref{prop:comparing-orbit-categories}, when defining the standard triangle structure on $\cd^b(\ch)/F$ in Theorem~\ref{thm:triangulated-orbit-category-Keller} we can sometimes replace a strict dg enhancement of $(\cd^b(\ch),F)$ involving $\cd^b(\ch)_{dg}$ by an equivalent strict dg enhancement involving $\cc_{dg}^b(\proj\ch)$. This is the case in Section~\ref{ss:twisted-root-categories}.
\end{remark}

Let $\ca$ be a dg $k$-category and $d\in\mathbb{N}$. Let $\cb$ be the direct sum of $d$ copies of $\ca$, namely, the objects of $\cb$ are $d$-tuples $(X_0,\ldots,X_{d-1})$ of objects of $\ca$, and for two objects $(X_0,\ldots,X_{d-1})$ and $(Y_0,\ldots,Y_{d-1})$ the morphism complex is 
\[
\Hom_\cb((X_0,\ldots,X_{d-1}),(Y_0,\ldots,Y_{d-1}))=\bigoplus_{i=0}^{d-1}\Hom_\ca(X_i,Y_i).
\]
Let $\Phi$ be a dg automorphism of $\ca$ and $\Psi$ be the dg automorphism of $\cb$ given by the matrix
\[
\left(\begin{array}{cccc}
& \id & & \\
& & \ddots &\\
& & & \id\\
\Phi & & &
\end{array}\right),
\]
where the rows and columns are labeled by $0,\ldots,d-1$.
Let $\iota\colon\ca\to\cb$ be the embedding of $\ca$ into $\cb$ as the $0$-th component. Then $\Psi^r\iota=\iota\Phi$ and $\iota$ induces a dg functor $\bar{\iota}\colon\ca/\Phi\to\cb/\Psi$. The proof of the following lemma is the same as that of \cite[Lemma 2.1]{Yang18}.

\begin{lemma}
\label{lem:comparing-orbit-categories}
$\bar{\iota}$ is a dg equivalence and the following diagram of dg functors is commutative
\[
\xymatrix{
\ca\ar[r]^\iota \ar[d]^\pi& \cb\ar[d]^\pi\\
\ca/\Phi\ar[r]^{\bar{\iota}} &\cb/\Psi.
}
\]
\end{lemma}

\section{Twisted root categories}
\label{s:twisted-root-categories}

In this section we introduce twisted root categories of Hom-finite and Ext-finite hereditary abelian categories and study in detail twisted root categories of two families of quivers of type $\mathbb{A}_\infty^\infty$.

\subsection{Twisted root categories}
\label{ss:twisted-root-categories}
Let $d\in\mathbb{Z}\backslash\{0\}$. Let $\ch$ be a small Hom-finite and Ext-finite hereditary abelian $k$-category and $\sigma$ be a $k$-linear automorphism of $\ch$. Then $\sigma$ induces a dg automorphism of $\cd^b(\ch)_{dg}$, the $H^0$ of which is exactly the derived functor of $\sigma$. We denote both functors still by $\sigma$. Then $\sigma[1]=[1]\sigma$. Consider the dg automorphism $\Phi=[d]\sigma$ of $\cd^b(\ch)_{dg}$. Then $H^0\Phi=\Sigma^d\sigma\colon \cd^b(\ch)\to\cd^b(\ch)$ satisfies the conditions 2) and 3) in Theorem~\ref{thm:triangulated-orbit-category-Keller}, so the orbit category $\cd^b(\ch)/\Sigma^d\sigma$ admits a triangle structure such that the projection functor $\pi\colon\cd^b(\ch)\to\cd^b(\ch)/\Sigma^d\sigma$ is a triangle functor. We call this triangulated category $\cd^b(\ch)/\Sigma^d\sigma$ the \emph{$d$-root category of $\ch$ twisted by $\sigma$}. It is clear that $\cd^b(\ch)/\Sigma^d\sigma=\cd^b(\ch)/\Sigma^{-d}\sigma^{-1}$. When $\sigma=\Id$, the $d$-root category $\cd^b(\ch)/\Sigma^d$ was introduced and studied in \cite{PengXiao97,PengXiao00,Zhang25} (with the name \emph{$d$-periodic derived category} in \cite{Zhang25}).

The shift functor of $\cd^b(\ch)/\Sigma^d\sigma$ is the functor induced from the shift functor $\Sigma$ of $\cd^b(\ch)$ with respect to the isomorphism $(-1)^d\id_{\Sigma^{d+1}\sigma}\colon\Sigma\Sigma^d\sigma\to\Sigma^d\sigma\Sigma$. This may cause confusion, because we also denote it by $\Sigma$ and there is a tendency to use $\id_{\Sigma^{d+1}\sigma}\colon\Sigma\Sigma^d\sigma\to\Sigma^d\sigma\Sigma$ instead of $(-1)^d\id_{\Sigma^{d+1}\sigma}$. The automorphism $\sigma$ of $\cd^b(\ch)$ with respect to the isomorphism $\id_{\Sigma^d\sigma^2}\colon\sigma\Sigma^d\sigma\to\Sigma^d\sigma\sigma$ induces an automorphism of $\cd^b(\ch)/\Sigma^d\sigma$, which we still denote by $\sigma$. Then we have $\Sigma^d(X)\cong\sigma^{-1}(X)$ in $\cd^b(\ch)/\Sigma^d\sigma$. A morphism in $\cd^b(\ch)$ is irreducible if and only if its image under $\pi$ is irreducible, and all irreducible morphisms in $\cd^b(\ch)/\Sigma^d\sigma$ are of this form modulo the square of the radical; a triangle in $\cd^b(\ch)$ is an AR-triangle if and only if its image under $\pi$ is an AR-triangle, and all AR-triangles in $\cd^b(\ch)/\Sigma^d\sigma$ are of this form up to isomorphism (see for example \cite[Section 2.3]{Yang18}). Therefore the Gabriel quiver (respectively, AR-quiver, when AR-triangles exist) of $\cd^b(\ch)/\Sigma^d\sigma$ is the orbit quiver of the Gabriel quiver (respectively, AR-quiver) of $\cd^b(\ch)$ by the action of $\Sigma^d\sigma$.

Assume that $\ch_0$ is an abelian subcategory of $\ch$ which is closed under $\sigma$. Then the derived functor $\iota\colon\cd^b(\ch_0)\to\cd^b(\ch)$ of the inclusion $\ch_0\hookrightarrow \ch$ is fully faithful with essential image $\cd^b_{\ch_0}(\ch)$, the full subcategory of $\cd^b(\ch)$ consisting of bounded complexes whose cohomologies belong to $\ch_0$. In this way we view $\cd^b(\ch_0)$ as a triangulated subcategory of $\cd^b(\ch)$. The inclusion $\ch_0\hookrightarrow \ch$ induces a dg functor $\cc^b_{dg}(\ch_0)\to\cc^b_{dg}(\ch)$, which takes acyclic complexes to acyclic complexes, and hence it further induces a dg functor $\tilde{\iota}\colon\cd^b(\ch_0)_{dg}\to\cd^b(\ch)_{dg}$ such that $H^0\tilde{\iota}=\iota$. The dg automorphism $\Phi$ is well-defined on $\cd^b(\ch_0)_{dg}$ and satisfies $\Phi\tilde{\iota}=\tilde{\iota}\Phi$. Then by Proposition~\ref{prop:comparing-orbit-categories}, the natural functor $\cd^b(\ch_0)/\Sigma^d\sigma\to\cd^b(\ch)/\Sigma^d\sigma$ induced by $\iota$ is a triangle functor. In this way, we consider $\cd^b(\ch_0)/\Sigma^d\sigma$ as a triangulated subcategory of $\cd^b(\ch)/\Sigma^d\sigma$. Therefore we have

\begin{lemma}
\label{lem:twisted-root-subcategory-of-subcategory}
There is a commutative diagram of triangle functors
\[
\xymatrix{
\cd^b(\ch)\ar[r]^\pi & \cd^b(\ch)/\Sigma^d\sigma\\
\cd^b(\ch_0)\ar[r]^\pi \ar@{^(->}[u] & \cd^b(\ch_0)/\Sigma^d\sigma. \ar@{^(->}[u]
}
\]
\end{lemma}

In the above process we also obtain a dg quasi-equivalence $\cd^b(\ch_0)_{dg}\to\cd^b_{\ch_0}(\ch)_{dg}$, the full dg subcategory of $\cd^b(\ch)_{dg}$ consisting of bounded complexes whose cohomologies belong to $\ch_0$.

Now assume that $\ch$ has enough projectives and let $\proj\ch$ be the full subcategory of $\ch$ of projective objects. By Remark~\ref{rem:triangle-structure-on-orbit-category-with-projectives}, when defining the standard triangle structure on $\cd^b(\ch)/\Sigma^d\sigma$ we can replace $\cd^b(\ch)_{dg}$ by $\cc_{dg}^b(\proj\ch)$. 
Assume further that $\ch_0$ is an abelian subcategory of $\ch$ which is closed under $\sigma$. Let $\cc_{dg,\ch_0}^b(\proj\ch)$ be the full dg subcategory of $\cc_{dg}^b(\proj\ch)$ consisting of bounded complexes whose cohomologies belong to $\ch_0$. Then the essential image of $H^0\cc_{dg,\ch_0}^b(\proj\ch)\hookrightarrow H^0\cc_{dg}^b(\proj\ch)=K^b(\proj\ch)\hookrightarrow \cd^b(\ch)$ is exactly $\cd^b_{\ch_0}(\ch)$. So there is a zigzag of quasi-equivalences
\[
\xymatrix{
\cc_{dg,\ch_0}^b(\proj\ch)\ar[r] & \cd^b_{\ch_0}(\ch)_{dg}\\
& \cd^b(\ch_0)_{dg}\ar[u]
}
\]
By Proposition~\ref{prop:comparing-orbit-categories}, when defining the standard triangle structure on $\cd^b(\ch_0)/\Sigma^d\sigma$ we can replace $\cd^b(\ch_0)_{dg}$ by $\cc_{dg,\ch_0}^b(\proj\ch)$.

\begin{lemma}
\label{lem:tilting-twisted-root-categories}
Let $\ch$ be a small Hom-finite and Ext-finite hereditary abelian $k$-category and $\sigma$ be a $k$-linear automorphism of $\ch$. Assume that $\cp\subset\ch$ is an idempotent-complete tilting subcategory such that there is a strongly locally finite quiver $Q$ with an automorphism which is also denoted by $\sigma$ together with a $k$-linear equivalence $\proj(Q)\to\cp$ which commutes with $\sigma$. Then
\begin{itemize}
\item[(a)] there is a triangle equivalence $\cd^b(\rep^+(Q))\to\cd^b(\ch)$, which commutes with $\sigma$;
\item[(b)] for $d\geq 1$ and $r\in\mathbb{Z}$, there is a triangle equivalence
\[
\cd^b(\rep^+(Q))/\Sigma^d\sigma^{-r}\longrightarrow\cd^b(\ch)/\Sigma^d\sigma^{-r}.
\]
\end{itemize}
\end{lemma}
\begin{proof}
There is a diagram of dg functors which commutes with $\sigma$
\[
\xymatrix{
\cc_{dg}^b(\proj(Q))\ar[r]\ar[d] &\cc_{dg}^b(\cp)\ar[r] & \cc_{dg}^b(\ch)\ar[d]\\
\cd^b(\rep^+(Q))_{dg} & & \cd^b(\ch)_{dg},
}
\] 
where the left vertical dg functor is a quasi-equivalence, and the first horizontal dg functor is the equivalence induced by the equivalence $\proj(Q)\to\cp$.
Taking $H^0$, we obtain a triangle functor $\cd^b(\rep^+(Q))\to\cd^b(\ch)$, which commutes with $\sigma$ and which is an equivalence by \cite[Lemma 1]{Beilinson78}. (b) is then obtained by Proposition~\ref{prop:comparing-orbit-categories}.
\end{proof}

We are particularly interested in twisted root categories of two families of quivers of type $\mathbb{A}_\infty^\infty$, which we study separately below.

\subsection{Special case: Quiver of type $\mathbb{A}_\infty^\infty$ with linear orientation} 
\label{ss:twisted-root-category-linear-orientation}

Let $Q^l$ be the quiver of type $\mathbb{A}_\infty^\infty$ with linear orientation:
\[\xymatrix{\cdots\ar[r]&i-1\ar[r]&i\ar[r]&i+1\ar[r]&\cdots}.\] 
Let $\sigma$ be the unique quiver-automorphism of $Q^l$ which takes $i$ to $i-1$. We denote the corresponding push-out functor also by $\sigma$. Below the translation quiver $\mathbb{Z}\mathbb{A}_\infty$ plays an important role. A useful feature is that on the underlying quiver of $\mathbb{Z}\mathbb{A}_\infty$ there is a unique structure of translation quiver.

\subsubsection{The Auslander--Reiten quiver}
The indecomposable representations of the category $\rep^b(Q^l)$ are precisely the $V_{a,b}$'s in Example~\ref{ex:representations-of-A-double-infinity}, where $a,b\in\mathbb{Z}$, $a\leq b$. It has Auslander--Reiten sequences and the Auslander--Reiten translation $\tau=\sigma^{-1}$ sends $V_{a,b}$ to $V_{a+1,b+1}$. Its AR-quiver $\cx$ is of shape $\mathbb{Z}\mathbb{A}_\infty$. The indecomposable representations of $\rep^+(Q^l)$ which are not in $\rep^b(Q^l)$ are precisely the indecomposable projective representations $P_a=V_{a,+\infty}$, $a\in\mathbb{Z}$, which form a connected component $\cy$ of the Gabriel quiver of $\rep^+(Q^l)$ of shape $(Q^l)^{op}$. Namely, the Gabriel quiver of $\rep^+(Q^l)$ consists of two components $\cy$ and $\cx$.
These results can be obtained by \cite[Lemma 4.5 (1) and Theorem 5.17 (1)]{BautistaLiuPaquette13} (see for example \cite[Proof of Theorem 4.16 (b)]{Yang18} for a detailed proof).
This implies, by \cite[Theorem 7.5]{BautistaLiuPaquette13}, that the bounded derived category $\cd^b(\rep^b(Q^l))$ has Auslander--Reiten triangles with $\tau=\sigma^{-1}$ and its AR quiver has $\mathbb{Z}$ connected components $\cx_i$ ($i\in\mathbb{Z}$) of shape $\mathbb{Z}\mathbb{A}_\infty$, and the Gabriel quiver of the bounded derived category $\cd^b(\rep^+(Q^l))$ has $\mathbb{Z}$ extra connected components $\cy_i$ ($i\in\mathbb{Z}$) of shape $(Q^l)^{op}$.

For $d\in\mathbb{N}$ and $r\in\mathbb{Z}$, consider $d$-root category of $\rep^+(Q^l)$ twisted by $\sigma^{-r}$:
\begin{align*}
\cc^+&=\cd^b(\rep^+(Q^l))/\Sigma^d\sigma^{-r},
\end{align*}
which, by Lemma~\ref{lem:twisted-root-subcategory-of-subcategory}, contains as a triangulated subcategory the $d$-root category of $\rep^b(Q^l)$ twisted by $\sigma^{-r}$:
\begin{align*}
\cc&=\cd^b(\rep^b(Q^l))/\Sigma^d\sigma^{-r}.
\end{align*}
Then $\cc$ has Auslander--Reiten triangles with $\tau=\sigma^{-1}$, and hence has a Serre functor $\mathbb{S}$, which is a $k$-linear endofunctor of $\cc$ such that there is a bifunctorial isomorphism $D\Hom_\cc(X,Y)\cong \Hom_\cc(Y,\mathbb{S}(X))$ for any $X,Y\in\cc$, where $D=\Hom_k(?,k)$ is the $k$-dual. It is unique up to a unique isomorphism and is determined  by the additive structure of $\cc$, and it extends to a triangle endofunctor of $\cc$. By \cite[Proposition I.2.3]{ReitenVandenBergh02}, $\mathbb{S}(X)\cong\Sigma\tau(X)$ for any $X\in\cc$. Recall that $\Sigma^d X\cong\sigma^r (X)$ for any $X\in\cc^+$.

\begin{lemma}
\label{lem:AR-quiver-of-twisted-root-category-linear-orientation} 
\begin{itemize}
\item[(a)] The AR-quiver of $\cc$ consists of $d$ copies $\cx_0,\ldots,\cx_{d-1}$ of  shape $\mathbb{Z}\mathbb{A}_\infty$. The shift functor $\Sigma$ identifies $\cx_{i}$ with $\cx_{i+1}$ (with the indices read modulo $d$).
\item[(b)] The Gabriel quiver of $\cc^+$ consists of $2d$ connected components: $\cy_0,\ldots,\cy_{d-1}$ of shape $(Q^l)^{op}$, and $\cx_0,\ldots,\cx_{d-1}$. The shift functor $\Sigma$ identifies $\cy_{i}$ with $\cy_{i+1}$ (with the indices read modulo $d$). 
\item[(c)] The Serre functor on $\cc$ satisfies $\mathbb{S}^r(X)\cong\Sigma^{-d+r}(X)$ and $\mathbb{S}^d(X)\cong\tau^{d-r}(X)$ for any $X\in\cc$.
\end{itemize}
\end{lemma}
\begin{proof}
(a) and (b) can be obtained from the description of the AR-quiver of $\cd^b(\rep^b(Q^l))$ and the Gabriel quiver of $\cd^b(\rep^+(Q^l))$. As for (c), we have 
\begin{align*}
\mathbb{S}^r(X)&\cong\Sigma^r\tau^r(X)\cong\Sigma^r\sigma^{-r}(X)\cong\Sigma^{-d+r}(X),\\
\mathbb{S}^d(X)&\cong\Sigma^d\tau^d(X)\cong\sigma^r\tau^d(X)\cong\tau^{d-r}(X).\qedhere
\end{align*}
\end{proof}

\subsubsection{Comparing twisted root categories}

\begin{proposition}
\label{prop:blowing-linear}
For $d,m\geq 1$ and $r\in\mathbb{Z}$, there is a commutative diagram of triangle functors
\begin{align*}
\xymatrix{
&\cd^b(\rep^+(Q^l))\ar[dl]\ar[dr]\\
\cd^b(\rep^+(Q^l))/\Sigma^{dm}\sigma^{-rm}\ar[rr] && \cd^b(\rep^+(Q^l))/\Sigma^d\sigma^{-r}\\
&\cd^b(\rep^b(Q^l))\ar@{-->}[dl]\ar@{-->}[dr]\ar@{^(-->}[uu]\\
\cd^b(\rep^b(Q^l))/\Sigma^{dm}\sigma^{-rm}\ar[rr]\ar@{^(->}[uu]&& \cd^b(\rep^b(Q^l))/\Sigma^d\sigma^{-r}\ar@{^(->}[uu]
}
\end{align*}
which realises $\cd^b(\rep^+(Q^l))/\Sigma^d\sigma^{-r}$ (respectively, $\cd^b(\rep^b(Q^l))/\Sigma^d\sigma^{-r}$) as an $m$-orbit category of $\cd^b(\rep^+(Q^l))/\Sigma^{dm}\sigma^{-rm}$ (respectively, $\cd^b(\rep^b(Q^l))/\Sigma^{dm}\sigma^{-rm}$).
\end{proposition}
\begin{proof}
This is obtained by taking $H^0$ of the commutative diagram \eqref{cd:finite-orbit-category} applied to $\cd^b(\rep^+(Q^l))_{dg}$ and $\cd^b(\rep^b(Q^l))_{dg}$ and $\Phi=[d]\sigma^{-r}$ in conjunction with Lemma~\ref{lem:twisted-root-subcategory-of-subcategory}.
\end{proof}

\subsubsection{Triangle equivalence vs additive equivalence}
\begin{proposition}
\label{prop:triangle-equivalence-vs-additive-equivalence-linear-orientation}
For $d,d'\in\mathbb{N}$ and $r,r'\in\mathbb{Z}$, the following conditions are equivalent:
\begin{itemize}
\item[(i)] there is a triangle equivalence $\cd^b(\rep^+(Q^l))/\Sigma^d\sigma^{-r}\to\cd^b(\rep^+(Q^l))/\Sigma^{d'}\sigma^{-r'}$,
\item[(ii)] there is a $k$-linear equivalence $\cd^b(\rep^+(Q^l))/\Sigma^d\sigma^{-r}\to\cd^b(\rep^+(Q^l))/\Sigma^{d'}\sigma^{-r'}$,
\item[(iii)] there is a triangle equivalence $\cd^b(\rep^b(Q^l))/\Sigma^d\sigma^{-r}\to\cd^b(\rep^b(Q^l))/\Sigma^{d'}\sigma^{-r'}$,
\item[(iv)] there is a $k$-linear equivalence $\cd^b(\rep^b(Q^l))/\Sigma^d\sigma^{-r}\to\cd^b(\rep^b(Q^l))/\Sigma^{d'}\sigma^{-r'}$,
\item[(v)] $(d,r)=(d',r')$.
\end{itemize}
\end{proposition}
\begin{proof}
The implications (v)$\Rightarrow$(i)$\Rightarrow$(ii) and (v)$\Rightarrow$(iii)$\Rightarrow$(iv) are obvious. It follows from Lemma~\ref{lem:AR-quiver-of-twisted-root-category-linear-orientation} (a) and (b) that (ii) implies (iv). Now assume that (iv) holds. Then the Gabriel quivers of $\cd^b(\rep^b(Q^l))/\Sigma^d\sigma^{-r}$ and $\cd^b(\rep^b(Q^l))/\Sigma^{d'}\sigma^{-r'}$ are isomorphic. By Lemma~\ref{lem:AR-quiver-of-twisted-root-category-linear-orientation} (a), this implies that $d=d'$. Then $r=r'$ follows by Lemma~\ref{lem:AR-quiver-of-twisted-root-category-linear-orientation} (c) from the fact that the Serre functor $\mathbb{S}$ is determined by the additive structure of the category and that in this case $\tau$ on objects is determined by the Gabriel quiver.
\end{proof}

\begin{remark}
Proposition~\ref{prop:triangle-equivalence-vs-additive-equivalence-linear-orientation} shows that if $\ct$ is a twisted root category $\cd^b(\rep^+(Q^l))/\Sigma^d\sigma^{-r}$ (respectively, $\cd^b(\rep^b(Q^l))/\Sigma^d\sigma^{-r}$), then its triangle structure is up to triangle equivalence determined by its additive structure. 
Moreover, for $d\in\mathbb{N}$ and $r,r'\in\mathbb{Z}$ with $r\neq r'$, the twisted root categories $\cd^b(\rep^+(Q^l))/\Sigma^d\sigma^{-r}$ and $\cd^b(\rep^+(Q^l))/\Sigma^{d}\sigma^{-r'}$ (respectively, $\cd^b(\rep^b(Q^l))/\Sigma^d\sigma^{-r}$ and $\cd^b(\rep^b(Q^l))/\Sigma^{d}\sigma^{-r'}$) are not $k$-linearly equivalent but they have isomorphic Gabriel quivers (respectively, AR-quivers). 
\end{remark}

\subsection{Special case: Quiver of type $\mathbb{A}_\infty^\infty$ with generalised zigzag orientation} 
\label{ss:twisted-root-category-zigzag-orientation}

Let $p\geq 0$ and $q\geq 1$ be integers. Let $Q_{p,q}$ be the quiver
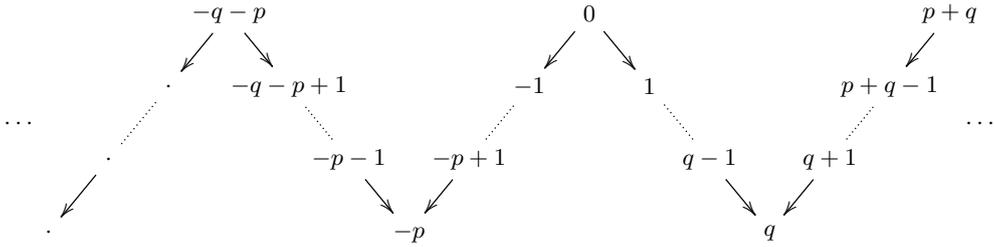
\[
{\scriptsize
\begin{xy} 0;<0.75pt,0pt>:<0pt,-0.55pt>::
(15,75) *+{\cdots}="",
(30,150) *+{\cdot}="0",
(60,100) *+{\cdot}="1",
(90,50) *+{\cdot}="2",
(120,0) *+{-q-p}="3",
(150,50) *+{-q-p+1}="4",
(180,100) *+{-p-1}="5",
(210,150) *+{-p}="6",
(240,100) *+{-p+1}="7",
(270,50) *+{-1}="8",
(300,0) *+{0}="9",
(330,50) *+{1}="10",
(360,100) *+{q-1}="11",
(390,150) *+{q}="12",
(420,100) *+{q+1}="13",
(450,50) *+{p+q-1}="14",
(480,0) *+{p+q}="15",
(495,75) *+{\cdots}="",
"1", {\ar "0"}, "2", {\ar@{.} "1"}, "3", {\ar "2"},
"3", {\ar "4"}, "4", {\ar@{.} "5"}, "5", {\ar "6"},
"7", {\ar "6"}, "8", {\ar@{.} "7"}, "9", {\ar "8"},
"9", {\ar "10"}, "10", {\ar@{.} "11"}, "11", {\ar "12"},
"13", {\ar "12"}, "14", {\ar@{.} "13"}, "15", {\ar "14"},
\end{xy}
}
\]
and $\sigma_{p,q}$ be the unique quiver-automorphism of $Q_{p,q}$ which takes $i$ to $i-q-p$. Note that $Q_{0,q}=Q^l$ and $\sigma_{0,q}=\sigma^q$.

In the rest of this subsection we assume $p\geq 1$.

\subsubsection{The Auslander--Reiten quiver}

By \cite[Theorems 4.6, 4.7 and 5.16]{BautistaLiuPaquette13}, $\rep^b(Q_{p,q})=\rep^+(Q_{p,q})$ has Auslander--Reiten sequences and its AR-quiver has four connected components: 
\begin{itemize}
\item[$\cdot$] the preprojective indecomposable representations form a component $\cp$ of shape $\mathbb{N}Q_{p,q}^{op}$,
\item[$\cdot$] the representations $\bigsqcup_{n\in\mathbb{Z}}\sigma_{p,q}^n\{V_{0,q},V_{-1,-1},\ldots,V_{-p+1,-p+1}\}$ and the iterated extensions of them form a regular component $\cx^1$ of shape $\mathbb{Z}A_\infty$,
\item[$\cdot$] the representations $\bigsqcup_{n\in\mathbb{Z}}\sigma_{p,q}^n\{V_{-p,0},V_{1,1},\ldots,V_{q-1,q-1}\}$ and the iterated extensions of them form a regular component $\cx^2$ of shape $\mathbb{Z}A_\infty$,
\item[$\cdot$] the preinjective indecomposable representations form a component $\ci$ of shape $(-\mathbb{N})Q_{p,q}^{op}$.
\end{itemize}
On the component $\cx^1$, we have $\tau^p=\sigma_{p,q}$; on the component $\cx^2$, we have $\tau^q=\sigma_{p,q}^{-1}$. 
By \cite[Theorem 7.10]{BautistaLiuPaquette13}, $\cd^b(\rep^b(Q_{p,q}))$ has Auslander--Reiten triangles and the AR quiver of $\cd^b(\rep^b(Q_{p,q}))$ has $3\mathbb{Z}$ connected components: $\cp_i$ ($i\in\mathbb{Z}$) of shape $\mathbb{Z}\mathbb{A}_\infty^\infty$ (which is glued from $\mathbb{N}Q_{p,q}^{op}$ and $(-\mathbb{N})Q_{p,q}^{op}$), $\cx^1_i$ ($i\in\mathbb{Z}$) of shape $\mathbb{Z}\mathbb{A}_\infty$ and $\cx^2_i$ ($i\in\mathbb{Z}$) of shape $\mathbb{Z}\mathbb{A}_\infty$.

For $d\in\mathbb{N}$ and $r\in\mathbb{Z}$, consider the $d$-root category of $\rep^b(Q_{p,q})$ twisted by $\sigma_{p,q}^{-r}$:
\begin{align*}
\cc&=\cd^b(\rep^b(Q_{p,q}))/\Sigma^d\sigma_{p,q}^{-r}.
\end{align*}
Then $\cc$ has Auslander--Reiten triangles, and hence has a Serre functor $\mathbb{S}$, which is isomorphic to $\Sigma\tau$ on objects, by \cite[Proposition I.2.3]{ReitenVandenBergh02}. Recall that $\Sigma^d X\cong\sigma^r (X)$ for any $X\in\cc$.

\begin{lemma}
\label{lem:AR-quiver-of-twisted-root-category-zigzag-orientation}
\begin{itemize}
\item[(a)] The AR-quiver of $\cc$ consists of $3d$ connected components: $\cp_0,\ldots,\cp_{d-1}$ of shape $\mathbb{Z}\mathbb{A}_\infty^\infty$, $\cx^1_0,\ldots,\cx^1_{d-1}$ of shape $\mathbb{Z}\mathbb{A}_\infty$ and $\cx^2_0,\ldots,\cx^2_{d-1}$ of shape $\mathbb{Z}\mathbb{A}_\infty$. The shift functor $\Sigma$ identifies $\cp_i$ with $\cp_{i+1}$, $\cx^1_i$ with $\cx^1_{i+1}$ and $\cx^2_i$ with $\cx^2_{i+1}$ (with the indices read modulo $d$).
\item[(b)] For $X\in\cx^1_i$ ($i=0,\ldots,d-1$), we have $\mathbb{S}^{pr}(X)\cong\Sigma^{pr+d}(X)$ and $\mathbb{S}^d(X)\cong\tau^{pr+d}(X)$; for $X\in\cx^2_i$ ($i=0,\ldots,d-1$), we have $\mathbb{S}^{qr}(X)\cong\Sigma^{qr-d}(X)$ and $\mathbb{S}^d(X)\cong\tau^{-qr+d}(X)$.
\end{itemize}
\end{lemma} \begin{proof}
(a) can be obtained from the description of the AR-quiver of $\cd^b(\rep^b(Q_{p,q}))$. As for (b), we have for $X\in\cx^1_i$ ($i=0,\ldots,d-1$)
\begin{align*}
\mathbb{S}^{pr}(X)&\cong\Sigma^{pr}\tau^{pr}(X)\cong\Sigma^{pr}\sigma_{p,q}^{r}(X)\cong\Sigma^{pr+d}(X),\\
\mathbb{S}^d(X)&\cong\Sigma^d\tau^d(X)\cong\sigma_{p,q}^r\tau^d(X)\cong\tau^{pr+d}(X),
\end{align*}
and for $X\in\cx^2_i$ ($i=0,\ldots,d-1$)
\begin{align*}
\mathbb{S}^{qr}(X)&\cong\Sigma^{qr}\tau^{qr}(X)\cong\Sigma^{qr}\sigma_{p,q}^{-r}(X)\cong\Sigma^{qr-d}(X),\\
\mathbb{S}^d(X)&\cong\Sigma^d\tau^d(X)\cong\sigma_{p,q}^r\tau^d(X)\cong\tau^{-qr+d}(X).\qedhere
\end{align*}
\end{proof}

\subsubsection{Comparing twisted root categories}
Let $p,q,d,r\geq 1$.

\begin{lemma}
\label{lem:turning-left-right}
\begin{itemize}
\item[(a)] The assignment $i\to -i$ extends uniquely to a quiver-isomorphism $\varphi\colon Q_{p,q}\to Q_{q,p}$ such that $\sigma_{q,p}^{-1}\varphi=\varphi\sigma_{p,q}$.
\item[(b)] The quiver-isomorphism $\varphi$ in (a) induces triangle equivalences
\begin{align*}
\cd^b(\rep^b(Q_{p,q}))&\longrightarrow\cd^b(\rep^b(Q_{q,p})),\\
\cd^b(\rep^b(Q_{p,q}))/\Sigma^d\sigma_{p,q}^{-r}&\longrightarrow\cd^b(\rep^b(Q_{q,p}))/\Sigma^d\sigma_{q,p}^r.
\end{align*}
\end{itemize}
\end{lemma}
\begin{proof}
(a) is clear. (b) follows by for example Lemma~\ref{lem:tilting-twisted-root-categories} with $\cp=\{P_i\mid i\in\mathbb{Z}\}$.
\end{proof}

\begin{proposition}
\label{prop:contract}
There are triangle equivalences
\begin{align*}
\cd^b(\rep^b(Q_{pr,qr}))/\Sigma^d\sigma_{pr,qr}^{-1}&{\longrightarrow} \cd^b(\rep^b(Q_{p,q}))/\Sigma^d\sigma_{p,q}^{-r},\\
\cd^b(\rep^b(Q_{qr,pr}))/\Sigma^d\sigma_{qr,pr}^{-1}&{\longrightarrow} \cd^b(\rep^b(Q_{p,q}))/\Sigma^d\sigma_{p,q}^{r}.
\end{align*}
\end{proposition}
\begin{proof}
Consider the additive closure  $\cp$ in $\rep^b(Q_{p,q})$ of 
\begin{align*}
\bigcup_{m\in\mathbb{Z}}\sigma_{p,q}^{rm}&(\{V_{-p,j+i(p+q)}\mid -p\leq j\leq -1,~ 0\leq i\leq r-1\}\\
&\cup\{V_{j+i(p+q),-p+r(p+q)}\mid -p\geq j\geq -p-q+1,~0\leq i\leq r-1\}).
\end{align*}
Then it is straightforward to check that $\cp$ is a tilting subcategory of $\rep^b(Q_{p,q})$, and that there is a $k$-linear equivalence $\proj(Q_{pr,qr})\to \cp$ which takes $\sigma_{pr,qr}$ to $\sigma_{p,q}^r$. So by Lemma~\ref{lem:tilting-twisted-root-categories}, we obtain the first desired triangle equivalence. The second triangle equivalence is obtained from the first one in conjunction with Lemma~\ref{lem:turning-left-right} (b).
\end{proof}

\begin{proposition}
\label{prop:blowing-zigzag}
$\cd^b(\rep^b(Q_{p,q}))/\Sigma^d\sigma_{p,q}^{-1}$ is an $r$-orbit category of  $\cd^b(\rep^b(Q_{p,q}))/\Sigma^{dr}\sigma_{p,q}^{-r}$, and there is a commutative diagram of triangle functors
\begin{align*}
\xymatrix{
\cd^b(\rep^b(Q_{p,q}))\ar[d]\ar[dr]\\
\cd^b(\rep^b(Q_{p,q}))/\Sigma^{dr}\sigma_{p,q}^{-r}\ar[r] & \cd^b(\rep^b(Q_{p,q}))/\Sigma^d\sigma_{p,q}^{-1}.
}
\end{align*}
\end{proposition}
\begin{proof}
This is obtained by taking $H^0$ of the commutative diagram \eqref{cd:finite-orbit-category} applied to $\ca=\cd^b(\rep^b(Q_{p,q}))_{dg}$ and $\Phi=[d]\sigma_{p,q}^{-1}$.
\end{proof}

\subsubsection{Triangle equivalence vs additive equivalence}
\begin{proposition}
\label{prop:triangle-equivalence-vs-additive-equivalence-zigzag-orientation}
For $p,q,d,p',q',d'\in\mathbb{N}$ and $r,r'\in\mathbb{Z}\backslash\{0\}$, the following conditions are equivalent:
\begin{itemize}
\item[(i)] there is a triangle equivalence $\cd^b(\rep^b(Q_{p,q}))/\Sigma^d\sigma_{p,q}^{-r}\to\cd^b(\rep^b(Q_{p',q'}))/\Sigma^{d'}\sigma_{p',q'}^{-r'}$,
\item[(ii)] there is a $k$-linear equivalence $\cd^b(\rep^b(Q_{p,q}))/\Sigma^d\sigma_{p,q}^{-r}\to\cd^b(\rep^b(Q_{p',q'}))/\Sigma^{d'}\sigma_{p',q'}^{-r'}$,
\item[(iii)] $(pr,qr,d)=(p'r',q'r',d')$ or $(pr,qr,d)=(-q'r',-p'r',d')$.
\end{itemize}
\end{proposition}
\begin{proof}
(i)$\Rightarrow$(ii): This is clear.

(iii)$\Rightarrow$(i): If $(pr,qr,d)=(p'r',q'r',d')$ or $(pr,qr,d)=(-q'r',-p'r',d')$, then by Proposition~\ref{prop:contract}, both $\cd^b(\rep^b(Q_{p,q}))/\Sigma^d\sigma_{p,q}^{-r}$ and $\cd^b(\rep^b(Q_{p',q'}))/\Sigma^{d'}\sigma_{p',q'}^{-r'}$ are triangle equivalent to $\cd^b(\rep^b(Q_{pr,qr}))/\Sigma^d\sigma_{pr,qr}^{-1}$ if $r\geq 1$ and to $\cd^b(\rep^b(Q_{-qr,-pr}))/\Sigma^d\sigma_{-qr,-pr}^{-1}$ if $r\leq -1$. Therefore they are triangle equivalent.

(ii)$\Rightarrow$(iii) Assume (ii). Then the Gabriel quivers of $\cd^b(\rep^b(Q_{p,q}))/\Sigma^d\sigma_{p,q}^{-r}$ and $\cd^b(\rep^b(Q_{p',q'}))/\Sigma^{d'}\sigma_{p',q'}^{-r'}$ are isomorphic. By Lemma~\ref{lem:AR-quiver-of-twisted-root-category-zigzag-orientation} (a) this implies that $d=d'$.  Then $(pr,qr)=(p'r',q'r')$ or $(-q'r',-p'r')$ follows from Lemma~\ref{lem:AR-quiver-of-twisted-root-category-zigzag-orientation} (b), from the fact that the Serre functor $\mathbb{S}$ is determined by the additive structure of the category and that in this case $\tau$ on objects is determined by the Gabriel quiver.
\end{proof}

\begin{remark}
Proposition~\ref{prop:triangle-equivalence-vs-additive-equivalence-zigzag-orientation} shows that if $\ct$ is a twisted root category $\cd^b(\rep^b(Q_{p,q}))/\Sigma^d\sigma_{p,q}^{-r}$, then its triangle structure is up to triangle equivalence determined by its additive structure. 
Moreover, for $d,p,q,p',q'\in\mathbb{N}$ and $r,r'\in\mathbb{Z}\backslash\{0\}$ with $(pr,qr)\neq (p'r',q'r')$ and $(pr,qr)\neq (-q'r',-p'r')$, the twisted root categories $\cd^b(\rep^b(Q_{p,q}))/\Sigma^d\sigma_{p,q}^{-r}$ and $\cd^b(\rep^b(Q_{p',q'}))/\Sigma^{d}\sigma_{p',q'}^{-r'}$ are not $k$-linearly equivalent but they have isomorphic AR-quivers. 
\end{remark}

\section{The derived category of a graded path algebra}
\label{s:derived-category-of-graded-path-algebra}

Let $Q$ be a finite connected graded quiver, and $kQ$ its graded path algebra. We view $kQ$ as a dg algebra with trivial differential and consider the perfect derived category $\per(kQ)$ of $kQ$ and the finite-dimensional derived category $\cd_{fd}(kQ)$. Let $\tilde{Q}$, $\mathrm{s}$ and $d$ be as defined in Section~\ref{ss:graded-quiver}. This section is devoted to describing $\per(kQ)$ and $\cd_{fd}(kQ)$ as certain twisted root categories associated to $\tilde{Q}$, in the case when $\per(kQ)$ is Hom-finite (and hence Krull--Schmidt).

We have $\cd_{fd}(kQ)\subseteq \per(kQ)$ by for example \cite[Section 3.14]{KellerYang11}, and the equality holds if and only if $Q$ has no oriented cycles.  This implies that the category $\cd_{fd}(kQ)$ is always Hom-finite. Indeed, it follows by d\'evissage that $\Hom_{\cd(kQ)}(M,X)$ is finite-dimensional for $M\in\per(kQ)$ and $X\in\cd_{fd}(kQ)$, because $\Hom_{\cd(kQ)}(kQ,X)=H^0(X)$ is finite-dimensional. Equivalent conditions for $\per(kQ)$ to be Hom-finite are given in Lemma~\ref{lem:finiteness-of-the-perfect-derived-category}.

\begin{theorem}
\label{thm:the-perfect-derived-category-of-graded-path-algebra}
Assume that $\per(kQ)$ is Hom-finite. 
If $d=0$, then there is a triangle equivalence 
\[
\xymatrix{
\cd^b(\rep(\widetilde{Q}))\ar[r]^(0.38){\simeq} & \per(kQ)=\cd_{fd}(kQ).
}
\]
If $d\geq 1$, then there is a commutative diagram of triangle functors with horizontal functors being triangle equivalences
\[
\xymatrix{
\cd^b(\rep^+(\widetilde{Q}))/\Sigma^d \mathrm{s}_*^{-1} \ar[r]^(0.62){\simeq}& \per(kQ) \\
\cd^b(\rep^b(\widetilde{Q}))/\Sigma^d \mathrm{s}_*^{-1}  \ar@{^{(}->}[u] \ar[r]^(0.62){\simeq} &   \cd_{fd}(kQ) \ar@{^{(}->}[u],
}
\]
where the triangulated categories on the left hand side are twisted root categories introduced in Section~\ref{s:twisted-root-categories} (\ie the triangle structures are the standard ones introduced in Section~\ref{ss:twisted-root-categories}).
\end{theorem}

\begin{remark}
\begin{itemize}
\item[(a)]
When $\per(kQ)$ is Hom-finite, the quiver $\tilde{Q}$ is strongly locally finite, see Lemma~\ref{lem:finiteness-of-the-perfect-derived-category}. Therefore the category $\rep^+(\widetilde{Q})$ is a Hom-finite and Ext-finite hereditary abelian $k$-category, by Proposition~\ref{prop:Hom-finiteness-of-finitely-presented-representations}.
\item[(b)] 
In the second statement, that the upper horizontal functor is a $k$-linear equivalence is implicitly contained in \cite[Section 4]{Yang18}. Also there it was claimed to be a triangle equivalence, but there the triangle structure on $\cd^b(\rep^+(\widetilde{Q}))/\Sigma^d \mathrm{s}_*^{-1}$ was the one induced via this equivalence.
\end{itemize}
\end{remark}

\begin{proof}
Let $i$ be the vertex which we fixed when we define $\tilde{Q}$.

First assume $d=0$. In this case, $Q$ has to be acyclic, and hence $kQ$ is finite-dimensional and $\per(kQ)=\cd_{fd}(kQ)$. Moreover, all walks from $i$ to a vertex $j$ are of the same degree, say $a_j$. Then $\bigoplus_{j\in Q_0} \Sigma^{a_j}e_jkQ$ is a tilting object of $\per(kQ)$ with endomorphism algebra $kQ^{un,op}$, where $Q^{un,op}$ is the opposite quiver of the underlying ungraded quiver of $Q$. By Lemma~\ref{lem:covering-quiver}, $\widetilde{Q}=Q^{un,op}$. Therefore there are triangle equivalences $\per(kQ)\simeq \per(k\widetilde{Q})\simeq K^b(\proj k\widetilde{Q})\cong\cd^b(\rep(\widetilde{Q}))$.

Next assume $d\geq 1$. By Section~\ref{ss:twisted-root-categories}, the orbit category $\cd^b(\grmod kQ)/\Sigma\langle-1\rangle$ admits a natural triangle structure and $\cd^b(\grmod_0 kQ)/\Sigma\langle-1\rangle$ is a triangulated subcategory. We claim that there is a triangle equivalence $\cd^b(\grmod kQ)/\Sigma\langle-1\rangle\to \per(kQ)$, which restricts to a triangle equivalence $\cd^b(\grmod_0 kQ)/\Sigma\langle-1\rangle\to \cd_{fd}(kQ)$. This is more or less proved in \cite[Theorem 5.1]{KalckYang18}. We only need a little more clarification. Indeed by \cite[Theorem 5.1]{KalckYang18}, the functor $\Tot\colon \cc_{dg}(\Grmod kQ)\to\cc_{dg}(kQ)$ of taking total complex induces a fully faithful dg functor $\overline{\Tot}\colon\cc_{dg}^b(\grproj kQ)/[1]\langle-1\rangle\to\cc_{dg}(kQ)$. Its essential image contains the dg modules $P_i=e_ikQ$, $i\in Q_0$, and is clearly contained in $\pretr(P_i\mid i\in Q_0)$. Thus we have the following commutative diagram of dg functors
\[
\xymatrix{
\cc_{dg}^b(\grproj kQ)/[1]\langle-1\rangle\ar[d]\ar[dr]^{\overline{\Tot}}&\\
\pretr(\cc_{dg}^b(\grproj kQ)/[1]\langle-1\rangle)\ar[r]& \pretr(P_i\mid i\in Q_0),
}
\]
where the horizontal functor is a dg equivalence. 
Taking $H^0$ we obtain a commutative diagram of fully faithful $k$-linear functors
\[
\xymatrix{
\cd^b(\grmod kQ)/\Sigma\langle-1\rangle\ar[d]\ar[dr]\ar[drr]&\\
H^0(\pretr(\cc_{dg}^b(\grproj kQ)/[1]\langle-1\rangle))\ar[r]& H^0\pretr(P_i\mid i\in Q_0)\ar[r] & \per(kQ),
}
\]
where the horizontal functors are triangle functors and the second one is an equivalence possibly up to direct summands.
The standard triangle structure on $\cd^b(\grmod kQ)/\Sigma\langle-1\rangle$ is obtained through the vertical functor, which is a $k$-linear equivalence by Theorem~\ref{thm:triangulated-orbit-category-Keller}, 
and by \cite[Theorem 5.1 (e)]{KalckYang18}, the right slanted functor is also a $k$-linear equivalence\footnote{In \cite[Theorem 5.1 (e)]{KalckYang18}, it is claimed to be a triangle equivalence. Actually there the triangle structure on the left side was the induced one via this equivalence.} which restricts to a $k$-linear equivalence $\cd^b(\grmod_0 kQ)/\Sigma\langle-1\rangle\to\cd_{fd}(kQ)$. This forces the horizontal functors to be triangle equivalences, it follows that the right slanted functor is a triangle equivalence and the restricted equivalence $\cd^b(\grmod_0 kQ)/\Sigma\langle-1\rangle\to\cd_{fd}(kQ)$ is also a triangle equivalence. 

Now we identify $R$ with $\coprod_d \tilde{Q}$ via the quiver-isomorphism $\varphi$ introduced in the end of Section~\ref{ss:graded-quiver}. Then there is a canonical isomorphism $\proj(R)=\proj(\coprod_d \tilde{Q})\cong\bigoplus_d \proj(\tilde{Q})$ of $k$-categories, which induces an isomorphism $\cc_{dg}^b(\proj(R))\cong\bigoplus_d\cc_{dg}^b(\proj(\tilde{Q}))$ of dg $k$-categories, restricting to an isomorphism $\cc_{dg,\rep^b}^b(\proj(R))\cong\bigoplus_d\cc_{dg,\rep^b}^b(\proj(\tilde{Q}))$.
The push-out dg functor $\sigma_*\colon \cc_{dg}^b(\proj(R))\to \cc_{dg}^b(\proj(R))$ corresponds to the dg automorphism of  $\bigoplus_d\cc_{dg}^b(\proj(\tilde{Q}))$ given by the matrix
\[
\left(\begin{array}{cccc}
& \id & & \\
& & \ddots &\\
& & & \id\\
\mathrm{s}_* & & &
\end{array}\right).
\] 
Thus applying Lemma~\ref{lem:comparing-orbit-categories} to $\ca=\cc_{dg}^b(\proj(\tilde{Q}))$, $\cb=\cc_{dg}^b(\proj(R))$, $\Psi=[1]\sigma_*^{-1}$ and $\Phi=[d] \mathrm{s}_*^{-1}$ we obtain a dg equivalence $
\cc_{dg}^b(\proj(\widetilde{Q}))/[d] \mathrm{s}_*^{-1}\to \cc_{dg}^b(\proj(R))/[1] \sigma_*^{-1}$, which restricts to a dg equivalence $
\cc_{dg,\rep^b}^b(\proj(\widetilde{Q}))/[d] \mathrm{s}_*^{-1}\to \cc_{dg,\rep^b}^b(\proj(R))/[1] \sigma_*^{-1}$. Moreover, a quasi-inverse of the restriction $\proj kQ\to\proj(R)$ of the $k$-linear equivalence $F$ in Lemma~\ref{lem:from-graded-modules-to-representations} induces a dg equivalence $\cc_{dg}^b(\proj(R))\to\cc_{dg}^b(\proj kQ)$, which takes $\sigma_*$ to $\langle1\rangle$ and which restricts to a dg equivalence $\cc_{dg,\rep^b}^b(\proj(R))\to\cc_{dg,\grmod_0}^b(\proj kQ)$. This further induces a dg equivalence $\cc_{dg}^b(\proj(R))/[1]\sigma_*^{-1}\to \cc_{dg}^b(\proj kQ)/[1]\langle-1\rangle$, which restricts to a dg equivalence $\cc_{dg,\rep^b}^b(\proj(R))/[1]\sigma_*^{-1}\to \cc_{dg,\grmod_0}^b(\proj kQ)/[1]\langle-1\rangle$. Composing the above dg equivalences and considering the pretriangulated hulls we obtain a commutative diagram of dg functors
\[
{\small
\begin{xy}0;<0.75pt,0pt>:<0pt,-0.75pt>::
(140,0) *+{\cc_{dg}^b(\proj(\widetilde{Q}))/[d] \mathrm{s}_*^{-1}}="0",
(420,0) *+{\cc_{dg}^b(\proj kQ)/[1] \langle-1\rangle}="1",
(0,50) *+{\cc_{dg,\rep^b}^b(\proj(\widetilde{Q}))/[d] \mathrm{s}_*^{-1}}="2",
(280,50) *+{\cc_{dg,\grmod_0}^b(\proj kQ)/[1] \langle-1\rangle}="3",
(140,130) *+{\pretr(\cc_{dg}^b(\proj(\widetilde{Q}))/[d] \mathrm{s}_*^{-1})}="4",
(420,130) *+{\pretr(\cc_{dg}^b(\proj kQ)/[1] \langle-1\rangle)}="5",
(0,180) *+{\pretr(\cc_{dg,\rep^b}^b(\proj(\widetilde{Q}))/[d] \mathrm{s}_*^{-1})}="6",
(280,180) *+{\pretr(\cc_{dg,\grmod_0}^b(\proj kQ)/[1] \langle-1\rangle)}="7",
"0", {\ar^(0.6)\simeq "1"}, "2", {\ar "0"}, "2", {\ar^(0.6)\simeq "3"}, "3", {\ar "1"},
"4", {\ar@{-->} "5"}, "6", {\ar@{-->} "4"}, "6", {\ar "7"}, "7", {\ar "5"},
"0", {\ar@{-->} "4"}, "1", {\ar "5"}, "2", {\ar "6"}, "3", {\ar "7"},
\end{xy}
}
\]
Taking $H^0$ we obtain a commutative diagram of $k$-linear functors
\[
{\small
\begin{xy}0;<0.75pt,0pt>:<0pt,-0.75pt>::
(140,0) *+{\cd^b(\rep^+(\widetilde{Q}))/\Sigma^d \mathrm{s}_*^{-1}}="0",
(420,0) *+{\cd^b(\grmod kQ)/\Sigma \langle-1\rangle}="1",
(0,50) *+{\cd^b(\rep^b(\widetilde{Q}))/\Sigma^d \mathrm{s}_*^{-1}}="2",
(280,50) *+{\cd^b(\grmod_0 kQ)/\Sigma \langle-1\rangle}="3",
(140,130) *+{H^0\pretr(\cc_{dg}^b(\proj(\widetilde{Q}))/[d] \mathrm{s}_*^{-1})}="4",
(420,130) *+{H^0\pretr(\cc_{dg}^b(\proj kQ)/[1] \langle-1\rangle)}="5",
(0,180) *+{H^0\pretr(\cc_{dg,\rep^b}^b(\proj(\widetilde{Q}))/[d] \mathrm{s}_*^{-1})}="6",
(280,180) *+{H^0\pretr(\cc_{dg,\grmod_0}^b(\proj kQ)/[1] \langle-1\rangle)}="7",
"0", {\ar^(0.6)\simeq "1"}, "2", {\ar "0"}, "2", {\ar^(0.6)\simeq "3"}, "3", {\ar "1"},
"4", {\ar@{-->} "5"}, "6", {\ar@{-->} "4"}, "6", {\ar "7"}, "7", {\ar "5"},
"0", {\ar@{-->} "4"}, "1", {\ar "5"}, "2", {\ar "6"}, "3", {\ar "7"},
\end{xy}
}
\]
The bottom square is a square of triangle functors.
By Theorem~\ref{thm:triangulated-orbit-category-Keller} the left two vertical functors are $k$-linear equivalences, which give the standard triangle structures on $\cd^b(\rep^+(\widetilde{Q}))/\Sigma^d \mathrm{s}_*^{-1}$ and on $\cd^b(\rep^b(\widetilde{Q}))/\Sigma^d \mathrm{s}_*^{-1}$. This shows that the top square is also a square of triangle functors, especially the two $k$-linear equivalences there are triangle equivalences. We obtain the desired triangle equivalences by composing these equivalences with those obtained in the first paragraph.
\end{proof}

\section{Graded gentle one-cycle algebras}
\label{s:graded-gentle-one-cycle-algebras}

In this section we study the perfect derived category and the finite-dimensional derived category of graded gentle one-cycle algebras (viewed as dg algebras with trivial diffential).
\medskip

A \emph{gentle algebra} is a finite-dimensional $k$-algebra of the form $kQ/(I)$, where $Q$ is a finite quiver and $I$ is a minimal set of relations, such that the following conditions hold:
\begin{itemize}
\item[-]
for each vertex of $Q$ there are at most two incoming arrows and at
most two outgoing arrows,
\item[-]
for each arrow $\beta$ of $Q$, both the number of arrows $\alpha$
with $t(\alpha)=s(\beta)$ and $\beta\alpha\notin I$ and the number
of arrows $\gamma$ with $s(\gamma)=t(\beta)$ and $\gamma\beta\notin
I$ are not greater than $1$,
\item[-]
for each arrow $\beta$ of $Q$, both the number of arrows $\alpha$
with $t(\alpha)=s(\beta)$ and $\beta\alpha\in I$ and the number of
arrows $\gamma$ with $s(\gamma)=t(\beta)$ and $\gamma\beta\in I$ are
not greater than $1$,
\item[-]
all relations in $I$ are paths of length $2$.
\end{itemize}
Gentle algebras are Gorenstein, due to \cite{GeissReiten05}. 
A \emph{graded gentle algebra} is a graded $k$-algebra whose underlying algebra is a gentle algebra.
A \emph{gentle one-cycle algebra} is a gentle algebra whose
underlying graph contains exactly one cycle. This cycle can be an oriented cycle, but if this is the case, then there has to be a relation on the cycle, due to the finite-dimensionality condition. A graded gentle one-cycle algebra $A$ belongs to one of the following two families:
\begin{itemize}
\item[(1)] $A$ is of finite global dimension. This happens if and only if the cycle is not an oriented cycle, or it is an oriented cycle  not with full relations. In this case $\per(A)=\cd_{fd}(A)$.
\item[(2)] $A$ is of infinite global dimension. This happens if and only if the cycle is an oriented cycle with full relations. In this case, $\per(A)\subset\cd_{fd}(A)$.
\end{itemize}

For a triple $(p,q,d)$ of integers such that $p\geq 0$ and $q\geq 1$, let $Q$ be the graded quiver
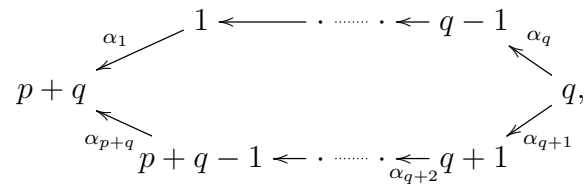
\[
\xymatrix@R=0.7pc@C=1pc{
&1\ar[dl]_{\alpha_{1}}&\cdot\ar[l]\ar@{.}[r]&\cdot&q-1\ar[l]&\\
p+q&&&&&q\ar[dl]^{\alpha_{q+1}}\ar[ul]_{\alpha_{q}},\\
&p+q-1\ar[ul]^{\alpha_{p+q}}&\cdot\ar[l]\ar@{.}[r]&\cdot&q+1\ar^(0.65){\alpha_{q+2}}[l]&
}
\]
where $\deg(\alpha_1)=d$ and the degrees of all the other arrows are $0$, and let $\Gamma(p,q,d)$ be its graded path algebra and $\hat{\Gamma}(p,q,d)$ be its complete graded path algebra. $\Gamma(p,q,d)$ and $\hat{\Gamma}(p,q,d)$ are different if and only if $p=0$ and $d=0$, \ie $Q$ is a cyclic quiver and the minimal oriented cycle is of degree $0$. 

The AG-invariant of a graded gentle algebra consists of a few pairs of integers. It was introduced in 
\cite{Avella-AlaminosGeiss08} for ungraded gentle algebras and in \cite{LekiliPolishchuk20} for graded gentle algebras.  It is a derived invariant for graded gentle algebras of finite global dimension, by \cite[Theorem 3.11]{LekiliPolishchuk20}, \cite[Theorem 1.4]{JinSchrollWang23} and \cite[Theorem C]{Opper25} (note that a graded gentle algebra shares the same geometric model and AG-invariant with its Koszul dual, which is a homologically smooth graded locally gentle algebra, by \cite[Section 9]{BoothGoodbodyOpper25}).

\begin{lemma}
\label{lem:AG-invariant-of-Gamma}
Let $p,q,d$ be integers such that $p\geq 0$ and $q\geq 1$.
\begin{itemize}
\item[(a)] For $p,q\geq 1$, the AG-invariant of $\Gamma(p,q,d)$ is $\{(p,p+d),(q,q-d)\}$.
\item[(b)] Let $\Gamma'(q,d)$ be the quotient of $\Gamma(0,q,q-d)$ modulo all paths of length $2$. Then the AG-invariant of $\Gamma'(q,d)$ is $\{(q,q-d),(0,d)\}$.
\end{itemize}
\end{lemma}
\begin{proof}
We refer to \cite[Section 3.1]{LekiliPolishchuk20} for unexplained terminologies and notation.

(a) There are $p+q$ forbidden threads: $\alpha_1,\alpha_2,\ldots,\alpha_{p+q}$ with grading
\[
|\alpha_1|_f=d,~~|\alpha_2|_f=\ldots=|\alpha_{p+q}|_f=0.
\]
Here we use $|\cdot|_f$ to denote the grading of a forbidden thread. There are $p+q$ permitted threads: $\alpha_1\cdots\alpha_q$, $\alpha_{p+q}\cdots\alpha_{q+1}$, $e_1,\ldots,e_{q-1}$, $e_{q+1},\ldots,e_{p+q-1}$ with grading
\begin{align*}
|\alpha_1\cdots\alpha_q|_p=-d,&~~|\alpha_{p+q}\cdots\alpha_{q+1}|_p=0,\\
|e_1|_p=\ldots=|e_{q-1}|_p=&~0=|e_{q+1}|_p=\ldots=|e_{p+q-1}|_p.
\end{align*}
Here we use $|\cdot|_p$ to denote the grading of a permitted thread.
Therefore there are two combinatorial boundary components, both of which are of type I:
\begin{align*}
b_1&=(\alpha_1\cdots\alpha_q)\alpha_{q+1}e_{q+1}\alpha_{q+2}\cdots e_{p+q-1}\alpha_{p+q},\\
b_2&=(\alpha_{p+q}\cdots\alpha_{q+1})\alpha_{q}e_{q-1}\alpha_{q-1}\cdots e_1\alpha_1,
\end{align*}
with $w(b_1)=-d$, $n(b_1)=p$, and $w(b_2)=d$, $n(b_2)=q$. Therefore the AG-invariant of $\Gamma(p,q,d)$ is $\{(p,p+d),(q,q-d)\}$.

(b) There are $2q$ forbidden threads: $\alpha_1\alpha_2\cdots\alpha_q$, $\alpha_2\cdots\alpha_q\alpha_1$, $\ldots$, $\alpha_q\alpha_1\cdots\alpha_{q-1}$, $e_1,e_2,\ldots,e_q$, with grading
\[
|\alpha_1\alpha_2\cdots\alpha_q|_f=\ldots=|\alpha_q\alpha_1\cdots\alpha_{q-1}|_f=1-d,~|e_1|_f=\ldots=|e_q|_f=1.
\]
There are $q$ permitted threads: $\alpha_1,\ldots,\alpha_q$ with grading
\[
|\alpha_1|_p=d-q,~|\alpha_2|_p=\ldots=|\alpha_q|_p=0.
\]
Therefore there are two combinatorial boundary components:
\begin{align*}
b_1&=\alpha_1e_1\alpha_2e_2\cdots\alpha_qe_q,\\
b_2&=\alpha_1\cdots\alpha_q,
\end{align*}
where $b_1$ is of type $I$ with $w(b_1)=d$ and $n(b_1)=q$, and $b_2$ is of type II' with $w(b_2)=-d$ and $n(b_2)=0$. Therefore the AG-invariant of $\Gamma'(q,d)$ is $\{(q,q-d),(0,d)\}$. 
\end{proof}

\begin{lemma}
\label{lem:canonical-form-of-graded-gentle-one-cycle-algebra}
Let $A$ be a graded gentle one-cycle algebra. 
\begin{itemize}
\item[(a)] If $A$ has finite global dimension, then its AG-invariant is $\{(p,p+d),(q,q-d)\}$ for some $p,q\geq 1$ and $d\in\mathbb{Z}$, and there is a triangle equivalence $\per(A)\to \per(\Gamma(p,q,d)))$.
\item[(b)] If $A$ has infinite global dimension, then its AG-invariant is $\{(q,q-d),(0,d)\}$ for some $q\geq 1$ and $d\in\mathbb{Z}$, and there are triangle equivalences $\cd_{fd}(A)\to\cd_{fd}(\Gamma'(q,d))\to\per(\hat{\Gamma}(0,q,d))$, which restrict to triangle equivalences $\per(A)\to\per(\Gamma'(q,d))\to \cd_{fd}(\hat{\Gamma}(0,q,d))$. 
\end{itemize}
\end{lemma}

\begin{proof}
(a)Assume that $A$ has finite global dimension. Then \cite[Theorem 1.1 (a)]{KalckYang18} states that there exists $p,q\geq 1$ and $d\in\mathbb{Z}$ such that $A$ and $\Gamma(p,q,d)$ are derived equivalent. Therefore there is a triangle equivalence $\per(A)\to \per(\Gamma(p,q,d))$, and moreover, the AG-invariant of $A$ is the same as the AG-invariant of $\Gamma(p,q,d)$, which is $\{(p,p+d),(q,q-d)\}$ by Lemma~\ref{lem:AG-invariant-of-Gamma} (a).

(b) Assume that $A$ has infinite global dimension. Then \cite[Theorem 1.1 (b)]{KalckYang18} says that there exists $q\geq 1$ and $d\in\mathbb{Z}$ such that $A$ and $\Gamma'(q,d)$ are derived equivalent. Therefore there is a triangle equivalence $\cd_{fd}(A)\to\cd_{fd}(\Gamma'(q,d))$, which restricts to a triangle equivalence $\per(A)\to \per(\Gamma'(q,d))$. The Koszul dual of $\Gamma'(q,d)$ is exactly $\hat{\Gamma}(0,q,d)$, and hence by \cite[Section 10.5 Lemma (the 'exterior` case)]{Keller94}, there is a triangle functor $\cd(\hat{\Gamma}(0,q,d))\to\cd(\Gamma'(q,d))$, which restricts to triangle equivalences $\per(\hat{\Gamma}(0,q,d))\to\cd_{fd}(\Gamma'(q,d))$ and $\cd_{fd}(\hat{\Gamma}(0,q,d))\to \thick(D\Gamma'(q,d))=\per(\Gamma'(q,d))$. Moreover, the AG-invariant of $A$ is the same as the AG-invariant of $\Gamma'(q,d)$, which is $\{(q,q-d),(0,d)\}$, by Lemma~\ref{lem:AG-invariant-of-Gamma} (b).
\end{proof}

\begin{remark}
\label{rem:AG-invariant}
\begin{itemize}
\item[(a)]
Note that $\Gamma(p,q,d)$ and $\Gamma(p',q',d')$ ($p,q,p,q'\geq 1$, $d,d'\in\mathbb{Z}$) have the same AG-invariants (up to permutation) if and only if $(p,q,d)=(p',q',d')$ or $(p,q,d)=(q',p',-d')$. Moreover, $\Gamma(p,q,d)$ and $\Gamma(q,p,-d)$ are graded equivalent, and hence derived equivalent. Similarly, $\Gamma'(q,d)$ and $\Gamma'(q',d')$ ($q,q'\geq 1$, $d,d'\in\mathbb{Z}$) has the same AG-invariants (up to permutation) if and only of $q=q'$ and $d=d'$. It follows then by \cite[Corollary 3.14 case (a)]{LekiliPolishchuk20}, \cite[Corollary 1.10]{LekiliPolishchuk20} and \cite[Theorem 1.4]{JinSchrollWang23}, \cite[Theorem C]{Opper25} that for graded gentle one-cycle algebras, the AG-invariant is a complete derived invariant, that is, two graded gentle one-cycle algebras of finite global dimension are derived equivalent if and only if their AG-invariants are the same (up to permutation).
\item[(b)]
A graded gentle one-cycle algebra $A$ is derived equivalent to an ungraded gentle one-cycle algebra if and only if the associated $(p,q,d)$ satisfies $-p\leq d\leq q$ when $A$ has finite global dimension, or the associated $(q,d)$ satisfies $1\leq d\leq q$ when $A$ has infinite global dimension. This follows from \cite[Theorem A]{BobinskiGeissSkowronski04} and \cite[Proposition 9.1]{KalckYang18}.
\end{itemize}
\end{remark}

For integers $p\geq 0$ and $q\geq 1$ we introduced in Section~\ref{ss:twisted-root-category-zigzag-orientation} a quiver $Q_{p,q}$ with quiver-automorphism $\sigma_{p,q}$. The next result describes the perfect derived category and the finite-dimensional derived category of $\Gamma(p,q,d)$ as twisted root categories.

\begin{lemma}
\label{lem:derived-category-of-Gamma-as-twisted-root-categories}
Let $(p,q,d)$ be a triple of integers such that  $p\geq 0,q\geq 1$ and $d\neq 0$. 
Then there is a triangle equivalence
\[
\xymatrix{
 \cd^b(\rep^+(Q_{p,q}))/\Sigma^d \sigma_{p,q}^{-1} \ar[r] & \per(\Gamma(p,q,d)),
}
\]
which restricts to a triangle equivalence
\[
\xymatrix{
 \cd^b(\rep^b(Q_{p,q}))/\Sigma^d \sigma_{p,q}^{-1} \ar[r] & \cd_{fd}(\Gamma(p,q,d)).
}
\]
\end{lemma}
\begin{proof}
We identify the quiver $\tilde{Q}$ with the connected component of $R$ containing the vertex $(p+q,0)$, which is
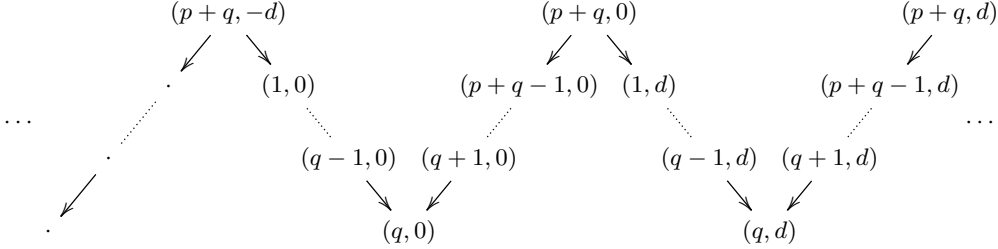
\[
{\scriptsize
\begin{xy} 0;<0.75pt,0pt>:<0pt,-0.55pt>::
(15,75) *+{\cdots}="",
(30,150) *+{\cdot}="0",
(60,100) *+{\cdot}="1",
(90,50) *+{\cdot}="2",
(120,0) *+{(p+q,-d)}="3",
(150,50) *+{(1,0)}="4",
(180,100) *+{(q-1,0)}="5",
(210,150) *+{(q,0)}="6",
(240,100) *+{(q+1,0)}="7",
(270,50) *+{(p+q-1,0)}="8",
(300,0) *+{(p+q,0)}="9",
(330,50) *+{(1,d)}="10",
(360,100) *+{(q-1,d)}="11",
(390,150) *+{(q,d)}="12",
(420,100) *+{(q+1,d)}="13",
(450,50) *+{(p+q-1,d)}="14",
(480,0) *+{(p+q,d)}="15",
(495,75) *+{\cdots}="",
"1", {\ar "0"}, "2", {\ar@{.} "1"}, "3", {\ar "2"},
"3", {\ar "4"}, "4", {\ar@{.} "5"}, "5", {\ar "6"},
"7", {\ar "6"}, "8", {\ar@{.} "7"}, "9", {\ar "8"},
"9", {\ar "10"}, "10", {\ar@{.} "11"}, "11", {\ar "12"},
"13", {\ar "12"}, "14", {\ar@{.} "13"}, "15", {\ar "14"},
\end{xy}
}
\]
and the quiver-automorphism $\mathrm{s}$ takes $(i,j)$ to $(i,j-|d|)$. Relabeling the vertices by $(i,j)\mapsto i+(\frac{j}{d}-1)(p+q)$, we obtain the quiver $Q_{p,q}$, and now the quiver-automorphism $\mathrm{s}$ takes $i$ to $i+p+q$ when $d\leq -1$ and to $i-p-q$ when $d\geq 1$.

By Theorem~\ref{thm:the-perfect-derived-category-of-graded-path-algebra}, there is a triangle equivalence
\[
\xymatrix{
 \cd^b(\rep^+(Q_{p,q}))/\Sigma^{|d|} \mathrm{s}_*^{-1}\ar[r] & \per(\Gamma(p,q,d)),
}
\]
which restricts a triangle equivalence
\[
\xymatrix{
 \cd^b(\rep^b(Q_{p,q}))/\Sigma^{|d|} \mathrm{s}_*^{-1}\ar[r] & \cd_{fd}(\Gamma(p,q,d)).
}
\]
If $d\geq 1$, then $|d|=d$ and $\sigma_{p,q}=\mathrm{s}$. If $d\leq -1$, then $|d|=-d$ and $\sigma_{p,q}=\mathrm{s}^{-1}$. In both cases, we obtain the desired result.
\end{proof}

\begin{remark}
\label{rem:d>=0}
\begin{itemize}
\item[(a)] It follows by Lemma~\ref{lem:turning-left-right} and Lemma~\ref{lem:derived-category-of-Gamma-as-twisted-root-categories} that there is a triangle equivalence $\per(\Gamma(p,q,d))\simeq\per(\Gamma(q,p,-d))$. This also follows from the fact that $\Gamma(p,q,d)$ and $\Gamma(q,p,-d)$ are graded equivalent. As a consequence, in Lemma~\ref{lem:canonical-form-of-graded-gentle-one-cycle-algebra} (a) we can always assume $d\geq 0$.
\item[(b)] Recall that $Q_{0,q}=Q^l$ and $\sigma_{0,q}=\sigma^q$. So
\begin{align*}
\cd^b(\rep^+(Q_{0,q}))/\Sigma^d \sigma_{0,q}^{-1}&=\cd^b(\rep^+(Q^l))/\Sigma^d\sigma^{-q}\\
&=\begin{cases} \cd^b(\rep^+(Q^l))/\Sigma^d\sigma^{-q} & \text{if $d\geq 1$}\\
\cd^b(\rep^+(Q^l))/\Sigma^{-d}\sigma^{q} & \text{if $d\leq -1$}
\end{cases},\\
\cd^b(\rep^b(Q_{0,q}))/\Sigma^d \sigma_{0,q}^{-1}&=\cd^b(\rep^b(Q^l))/\Sigma^d\sigma^{-q}\\
&=\begin{cases} \cd^b(\rep^b(Q^l))/\Sigma^d\sigma^{-q} & \text{if $d\geq 1$}\\
\cd^b(\rep^b(Q^l))/\Sigma^{-d}\sigma^{q} & \text{if $d\leq -1$}
\end{cases}.
\end{align*} 
\end{itemize}
\end{remark}

\begin{theorem}
\label{thm:derived-category-as-twisted-root-category}
Let $A$ be a graded gentle one-cycle algebra.
\begin{itemize}
\item[(a)] If the AG-invariant of $A$ is $\{(p,p),(q,q)\}$ for some $p,q\geq 1$, then there is a triangle equivalence $\per(A)\to\cd^b(\rep^b(Q))$, where $Q$ is a/any finite quiver of type $\tilde{\mathbb{A}}_{p,q}$. 

\item[(b)] If the AG-invariant of $A$ is $\{(p,p+d),(q,q-d)\}$ for some $p,q,d\geq 1$, then there is a triangle equivalence $\per(A)\to\cd^b(\rep^b(Q_{p,q}))/\Sigma^d\sigma_{p,q}^{-1}$.

\item[(c)]If the AG-invariant of $A$ is $\{(q,q),(0,0)\}$ for some $q\geq 1$, then there is a triangle equivalence $\cd_{fd}(A)\to\cd^b(\mod\hat{kQ})$ which restricts to a triangle equivalence $\per(A)\to\cd^b(\mod_0\hat{kQ})\simeq\cd^b(\rep_{\rm nil}^b(Q))$, where $Q$ is the cyclic quiver with $q$ vertices. Here $\mod\hat{kQ}$ is the category of finitely presented modules over $\hat{kQ}$, $\mod_0\hat{kQ}$ is the category of finite-dimensional modules over $\hat{kQ}$ and $\rep_{\rm nil}^b(Q)$ is the category of finite-dimensional nilpotent representations of $Q$. 

\item[(d)] If the AG-invariant of $A$ is $\{(q,q-d),(0,d)\}$ for some $q\geq 1$ and $d\neq 0$, then there is a triangle equivalence $\cd_{fd}(A)\to\cd^b(\rep^+(Q^l))/\Sigma^{|d|}\sigma^{-\mathrm{sgn}(d)q}$ which restricts to a triangle equivalence $\per(A)\to\cd^b(\rep^b(Q^l))/\Sigma^{|d|}\sigma^{-\mathrm{sgn}(d)q}$, where $\mathrm{sgn}(d)$ is the sign of $d$.
\end{itemize}
\end{theorem}
\begin{proof}
(a) Assume that the AG-invariant of $A$ is $\{(p,p),(q,q)\}$. Then it follows by Lemma~\ref{lem:canonical-form-of-graded-gentle-one-cycle-algebra} (a) that there is a triangle equivalence $\per(A)\to \per(\Gamma(p,q,0)))$. Note that $\Gamma(p,q,0)$ is the path algebra of a quiver of type $\tilde{\mathbb{A}}_{p,q}$, and hence there is a triangle equivalence $\per(A)\to\cd^b(\rep^b(Q))$, where $Q$ is any finite quiver of type $\tilde{\mathbb{A}}_{p,q}$. 

(b) Assume that the AG-invariant of $A$ is $\{(p,p+d),(q,q-d)\}$ with $d\geq 1$. Then it follows by Lemma~\ref{lem:canonical-form-of-graded-gentle-one-cycle-algebra} (a) that there is a triangle equivalence $\per(A)\to \per(\Gamma(p,q,d)))$, and the latter category is triangle equivalent to $\cd^b(\rep^b(Q_{p,q}))/\Sigma^d\sigma_{p,q}^{-1}$ by Lemma~\ref{lem:derived-category-of-Gamma-as-twisted-root-categories}. Composing the two triangle equivalences we obtain a triangle equivalence $\per(A)\to\cd^b(\rep^b(Q_{p,q}))/\Sigma^d\sigma_{p,q}^{-1}$.

(c) Assume that the AG-invariant of $A$ is $\{(q,q),(0,0)\}$ for some $q\geq 1$. Then it follows by Lemma~\ref{lem:canonical-form-of-graded-gentle-one-cycle-algebra} (b) that there is a triangle equivalence $\cd_{fd}(A)\to\per(\hat{\Gamma}(0,q,0))$ which restricts to a triangle equivalence $\per(A)\to\cd_{fd}(\hat{\Gamma}(0,q,0))$. Here $\hat{\Gamma}(0,q,0)$ is the complete path algebra $\hat{kQ}$ of the cyclic quiver $Q$ with $q$ vertices, and $\per(\hat{\Gamma}(0,q,0))=\cd^b(\mod\hat{kQ})$ and $\cd_{fd}(\hat{\Gamma}(0,q,0))=\cd^b(\mod_0\hat{\Gamma}(0,q,0))=\cd^b(\mod_0\hat{kQ})$, where $\mod_0\hat{kQ}$ is known to be isomorphic to $\rep^b_{\rm nil}(Q)$.

(d) Assume that the AG-invariant of $A$ is $\{(q,q-d),(0,d)\}$ for some $q\geq 1$ and $d\neq 0$, then it follows by Lemma~\ref{lem:canonical-form-of-graded-gentle-one-cycle-algebra} (b) that there is a triangle equivalence $\cd_{fd}(A)\to\per(\hat{\Gamma}(0,q,d))$ which restricts to a triangle equivalence $\per(A)\to\cd_{fd}(\hat{\Gamma}(0,q,d))$. Here $\hat{\Gamma}(0,q,d)=\Gamma(0,q,d)$ and it follows from Lemma~\ref{lem:derived-category-of-Gamma-as-twisted-root-categories} that there is a triangle equivalence $\cd_{fd}(A)\to\cd^b(\rep^+(Q_{0,q}))/\Sigma^d\sigma_{0,q}^{-1}$ which restricts to a triangle equivalence $\per(A)\to\cd^b(\rep^b(Q_{0,q}))/\Sigma^d\sigma_{0,q}^{-1}$. This finishes the proof in view of Remark~\ref{rem:d>=0} (b).
\end{proof}

In fact we have the following characterisation of perfect derived categories of graded gentle one-cycle algebras.

\begin{theorem}
\label{thm:perfect-derived-category-as-twisted-root-category}
Let $\ct$ be a triangulated $k$-category. The following conditions are equivalent:
\begin{itemize}
\item[(i)] $\ct$ is the perfect derived category of a graded gentle one-cycle algebra,
\item[(ii)] $\ct$ is triangle equivalent to one triangulated category in the following four families:
\begin{itemize}
\item[(1)] $\cd^b(\rep^b(Q^l))/\Sigma^d\sigma^{-r}$, where $d\geq 1$ and $r\in\mathbb{Z}\backslash\{0\}$, 
\item[(2)] $\cd^b(\rep^b(Q_{p,q}))/\Sigma^d\sigma_{p,q}^{-r}$, where $p,q,d\geq 1$ and $r\in\mathbb{Z}\backslash\{0\}$, 
\item[(3)] $\cd^b(\mod_0\hat{kQ})\simeq\cd^b(\rep_{\rm nil}^b(Q))$, where $Q$ is a cyclic quiver,
\item[(4)] $\cd^b(\rep^b(Q))$, where $Q$ is a finite quiver of type $\tilde{\mathbb{A}}_{p,q}$, $p,q\geq 1$.
\end{itemize}
\end{itemize}
\end{theorem}
\begin{proof}
(i)$\Rightarrow$(ii) This follows from Theorem~\ref{thm:derived-category-as-twisted-root-category}.
 
 (ii)$\Rightarrow$(i) By Remark~\ref{rem:d>=0} (b),  Lemma~\ref{lem:derived-category-of-Gamma-as-twisted-root-categories} and Lemma~\ref{lem:canonical-form-of-graded-gentle-one-cycle-algebra} (b), a twisted root category $\cd^b(\rep^b(Q^l))/\Sigma^d\sigma^{-r}$ in the family (1) is triangle equivalent to $\per(\Gamma'(|r|,(-1)^{\mathrm{sgn}(r)}d))$, and a category in the family (3) is triangle equivalent to $\per(\Gamma'(q,0))$, where $q$ is the number of vertices of $Q$, in view of the proof of Theorem~\ref{thm:derived-category-as-twisted-root-category}. By Proposition~\ref{prop:contract} and Lemma~\ref{lem:derived-category-of-Gamma-as-twisted-root-categories}, a twisted root category $\cd^b(\rep^b(Q_{p,q}))/\Sigma^d\sigma_{p,q}^{-r}$ in the family (2) is triangle equivalent to $\per(\Gamma(pr,qr,d))$ if $r\geq 1$ and to $\per(\Gamma(-qr,-pr,d))$ if $r\geq -1$. For the family (4), it is enough to know that the path algebra of any quiver of type $\tilde{\mathbb{A}}_{p,q}$ is a gentle algebra.
\end{proof}

We see from the proof of Theorem~\ref{thm:perfect-derived-category-as-twisted-root-category} that for a triangulated category $\ct$ in (ii), the AG-invariant of the corresponding graded gentle one-cycle algebra can be computed from the integers $p,q,d,r$.

\smallskip
The following is a `finite-dimensional' version of Theorem~\ref{thm:perfect-derived-category-as-twisted-root-category}. The proof is similar, and we omit it.

\begin{theorem}
\label{thm:finite-dimensional-derived-category-as-twisted-root-category}
Let $\ct$ be a triangulated $k$-category. The following conditions are equivalent:
\begin{itemize}
\item[(i)] $\ct$ is the finite-dimensional derived category of a graded gentle one-cycle algebra,
\item[(ii)] $\ct$ is triangle equivalent to one triangulated category in the following four families:
\begin{itemize}
\item[(1)] $\cd^b(\rep^+(Q^l))/\Sigma^d\sigma^{-r}$, where $d\geq 1$ and $r\in\mathbb{Z}\backslash\{0\}$, 
\item[(2)] $\cd^b(\rep^b(Q_{p,q}))/\Sigma^d\sigma_{p,q}^{-r}$, where $d,p,q\geq 1$ and $r\in\mathbb{Z}\backslash\{0\}$, 
\item[(3)] $\cd^b(\mod\hat{kQ})$, where $Q$ is a cyclic quiver,
\item[(4)] $\cd^b(\rep^b(Q))$, where $Q$ is a finite quiver of type $\tilde{\mathbb{A}}_{p,q}$, $p,q\geq 1$.
\end{itemize}
\end{itemize}
\end{theorem}

The following result shows that for if $\ct$ is the perfect derived category (respectively, the finite-dimensional derived category) of a graded gentle one-cycle algebra, its triangle structure is up to triangle equivalence determined by the additive category. Actually, for a triangulated category $\ct$ in Theorem~\ref{thm:perfect-derived-category-as-twisted-root-category} (ii), the integers $p,q,d,r$ can be computed from the additive structure of $\ct$. 

\begin{corollary}
\label{cor:additive=>triangle}
Let $A$ and $A'$ be two graded gentle one-cycle algebras. Then the following conditions are equivalent:
\begin{itemize}
\item[(i)] there is a triangle equivalence $\per(A)\to\per(A')$,
\item[(ii)] there is a $k$-linear equivalence $\per(A)\to\per(A')$,
\item[(iii)] there is a triangle equivalence $\cd_{fd}(A)\to\cd_{fd}(A')$,
\item[(iv)] there is a $k$-linear equivalence $\cd_{fd}(A)\to\cd_{fd}(A')$.
\end{itemize}
\end{corollary}
\begin{proof}
This follows from Theorem~\ref{thm:derived-category-as-twisted-root-category}, Propositions~\ref{prop:triangle-equivalence-vs-additive-equivalence-linear-orientation}~and~\ref{prop:triangle-equivalence-vs-additive-equivalence-zigzag-orientation} and Lemma~\ref{lem:the-derived-category-of-a-cycle} below.
\end{proof}

\begin{lemma}
\label{lem:the-derived-category-of-a-cycle}
Let $Q$ and $Q'$ be two cyclic quivers with $q$ and $q'$ vertices, respectively. Then the following conditions are equivalent:
\begin{itemize}
\item[(i)] there is a triangle equivalence $\cd^b(\mod\hat{kQ})\to\cd^b(\mod\hat{kQ'})$,
\item[(ii)] there is a $k$-linear equivalence $\cd^b(\mod\hat{kQ})\to\cd^b(\mod\hat{kQ'})$,
\item[(iii)] there is a triangle equivalence $\cd^b(\mod_0\hat{kQ})\to\cd^b(\mod_0\hat{kQ'})$,
\item[(iv)] there is a $k$-linear equivalence $\cd^b(\mod_0\hat{kQ})\to\cd^b(\mod_0\hat{kQ'})$,
\item[(v)] $q=q'$.
\end{itemize}
\end{lemma}
\begin{proof}
The implications (v)$\Rightarrow$(i)$\Rightarrow$(ii) and (v)$\Rightarrow$(iii)$\Rightarrow$(iv) are clear.

(ii)$\Rightarrow$(v) This is because the Gabriel quiver of $\mod\hat{kQ}$ consists of two connected components: one is of shape $Q$ and the other is a tube of rank $q$, and Gabriel quiver of $\cd^b(\mod\hat{kQ})$ consists of $\mathbb{Z}$ connected components of shape $Q$ and $\mathbb{Z}$ connected components which are tubes of rank $q$.

(iv)$\Rightarrow$(v) This is because the Gabriel quiver of $\mod_0\hat{kQ}$ is a tube of rank $q$, and the Gabriel quiver of $\cd^b(\mod_0\hat{kQ})$ consists of $\mathbb{Z}$ connected components which are tubes of rank $q$.
\end{proof}

The following corollary is a consequence of Theorem~\ref{thm:derived-category-as-twisted-root-category} and Propositions~\ref{prop:blowing-linear} and \ref{prop:blowing-zigzag}. When $p\geq 1$, the first row and the second row of the following diagram coincide.

\begin{corollary}
\label{cor:comparing-graded-gentle-one-cycle-algebras}
Let $p\geq 0$, $q,r\geq 1$ and $d\neq 0$ be integers. If $A$ and $A'$ are graded gentle one-cycle algebras with AG-invariants $\{(pr,pr+dr),(qr,qr-dr)\}$ and $\{(p,p+d),(q,q-d)\}$, respectively, then there is a commutative diagram of triangle functors
\[
\xymatrix{
\cd_{fd}(A) \ar[r] & \cd_{fd}(A')\\
\per(A) \ar[r]\ar@{^(->}[u] & \per(A')\ar@{^(->}[u]
}
\]
which realises $\cd_{fd}(A')$ (respectively, $\per(A')$) as an $r$-orbit category of $\cd_{fd}(A)$ (respectively, $\per(A)$).
\end{corollary}

\def\cprime{$'$}
\providecommand{\bysame}{\leavevmode\hbox to3em{\hrulefill}\thinspace}
\providecommand{\MR}{\relax\ifhmode\unskip\space\fi MR }
\providecommand{\MRhref}[2]{%
  \href{http://www.ams.org/mathscinet-getitem?mr=#1}{#2}
}
\providecommand{\href}[2]{#2}

\end{document}